\documentclass{amsart}
\usepackage{amssymb,delimset,booktabs,keytheorems,graphicx,tabularx,caption,floatrow}
\usepackage{mathtools}

\usepackage{enumitem}
\setlist{itemsep=4pt, topsep=0pt, leftmargin=17pt}

\usepackage{xcolor}

\definecolor{BIT}{cmyk}{1, 0, 1, 0}

\usepackage[numbers, sort&compress, nonamebreak, merge, elide, longnamesfirst]{natbib}

\usepackage{tikz}
\usetikzlibrary{decorations.pathreplacing, calc, positioning, patterns}
\tikzset{
edge/.style={semithick}}

\usepackage[pagebackref,
colorlinks=true,
linkcolor=blue,
citecolor=blue,
urlcolor=blue,
linktoc=all,
bookmarks=true]{hyperref}

\usepackage[capitalize,sort&compress,nameinlink]{cleveref}
\crefname{figure}{Fig.}{Figs.}
\Crefname{figure}{Fig.}{Figs.}
\crefname{table}{Tab.}{Tabs.}
\Crefname{table}{Tab.}{Tabs.}
\crefname{appendix}{Appendix}{Appendices}
\Crefname{appendix}{Appendix}{Appendices}
\crefformat{appendix}{#2Appendix~#1#3}
\Crefformat{appendix}{#2Appendix~#1#3}
\crefname{itm}{Condition}{Conditions}
\creflabelformat{itm}{~\upshape(#2#1#3)}
\crefname{def}{Def.}{Defs.}
\Crefname{def}{Definition}{Definitions}
\creflabelformat{def}{~\upshape(#2#1#3)}
\crefname{ineq}{Ineq.}{Ineqs.}
\Crefname{ineq}{Inequality}{Inequalities}
\creflabelformat{ineq}{~\upshape(#2#1#3)}
\crefname{conjecture}{Conjecture}{Conjectures}
\crefformat{section}{#2\S#1#3}
\Crefformat{section}{#2\S#1#3}
\crefrangeformat{section}{#3\S#1#4--#5\S#2#6}
\crefmultiformat{section}
{#2\S#1#3}
{ and~#2\S#1#3}
{, #2\S#1#3}
{, and~#2\S#1#3}

\AddToHook{env/conjecture/begin}{\crefalias{theorem}{conjecture}}
\AddToHook{env/corollary/begin}{\crefalias{theorem}{corollary}}
\AddToHook{env/lemma/begin}{\crefalias{theorem}{lemma}}
\AddToHook{env/proposition/begin}{\crefalias{theorem}{proposition}}

\makeatletter
\newcommand\creflabel[2][\@currentcounter]{%
 \crefalias{\@currentcounter}{#1}\label{#2}}
\makeatother

\newtheorem{theorem}{Theorem}[section]
\newkeytheorem{corollary}[name=Corollary, sibling=theorem]
\newkeytheorem{conjecture}[name=Conjecture, sibling=theorem]
\newkeytheorem{lemma}[name=Lemma, sibling=theorem]
\newkeytheorem{proposition}[name=Proposition, sibling=theorem]
\newkeytheorem{problem}[name=Problem, sibling=theorem]

\theoremstyle{definition}

\numberwithin{equation}{section}
\numberwithin{table}{section}
\numberwithin{figure}{section}

\allowbreak
\allowdisplaybreaks

\DeclareMathOperator{\fp}{fp}
\DeclareMathOperator{\last}{last}
\DeclareMathOperator{\wt}{wt}

\newcommand{\Cge}[1]{\mathcal{C}_{#1}^{\ge 2}}
\newcommand{\Cle}[1]{\mathcal{C}_{#1}^{\le 2}}

\title[A noncommutative approach to Schur positivity]
{Schur positivity of the spiders~$S(a,2,1)$ and~$S(a,4,1)$ via noncommutative symmetric functions}

\author[J.-Y.~Thibon]{Jean-Yves Thibon}
\address[Jean-Yves Thibon]{
LIGM,
Université Gustave-Eiffel,
CNRS,
ENPC,
ESIEE-Paris\\
5 Boulevard Descartes\\
Champs-sur-Marne 77454 Marne-la-Vallée cedex 2,
France}
\email{jean-yves.thibon@univ-eiffel.fr}
\thanks{Thibon is supported by ANR Grant CARPLO (ANR-20-CE40-0007).}

\author[D.G.L.~Wang]{David G.L. Wang$^*$}
\address[David G.L. Wang]{
School of Mathematics and Statistics \& MIIT Key Laboratory of Mathematical Theory and Computation in Information Security\\
Beijing Institute of Technology\\
Beijing 102400, P. R. China}
\email{glw@bit.edu.cn}
\thanks{$^*$Wang is the corresponding author, and is supported by the NSFC (Grant No.\ 12171034).}

\date{}

\hypersetup{
pdftitle={Schur positivity of the spiders S(a,2,1) and S(a,4,1) via noncommutative symmetric functions},
pdfauthor={Jean-Yves Thibon and David G.L. Wang},
pdfkeywords={chromatic symmetric functions, e-positivity, Littlewood--Richardson rule, noncommutative symmetric functions, Schur positivity}
}

\begin{document}

\begin{abstract}
We prove that the spider graphs~$S(a,2,1)$ and~$S(a,4,1)$ are Schur positive for all integers~$a\ge1$.
Together with the known~$e$-positivity results,
this completes the~$e$- and Schur-positivity classification of both families.
Our approach uses noncommutative symmetric functions,
including a particularly simple ribbon expansion for the path lift with coefficients given by powers of two.
We give a new proof of the Shareshian--Wachs path formula at~$t=1$ and construct corresponding lifts for spiders.
The Littlewood--Richardson rule converts their ribbon expansions into a general Schur-coefficient formula
in terms of weighted Yamanouchi words.
Mass-preserving multi-injections and a reduction to finitely many inequalities in degree~$10$
then prove the required positivity;
exact computer verification of these inequalities completes the proof in full generality.
\end{abstract}

\maketitle

\section{Introduction}

\citet{Sta95} introduced the \emph{chromatic symmetric function} of a graph~$G$ as
\[
X_G(x_1,x_2,\dots)
=\sum_{\text{proper colorings $c$ of $G$}}
\prod_{v\in V(G)}x_{c(v)}.
\]
This is a symmetric-function generalization of Birkhoff's chromatic polynomial.
A graph~$G$ is~\emph{$e$-positive} if the expansion of~$X_G$ in the elementary basis has nonnegative coefficients,
and \emph{Schur positive} if its expansion in the Schur basis has nonnegative coefficients.
By the Pieri rule, any $e$-positive symmetric function is Schur positive.

Two conjectures have shaped the study of these positivity properties.
\citet{SS93} conjectured that the incomparability graphs of~$(3+1)$-free posets are~$e$-positive.
An equivalent statement of the Stanley--Stembridge conjecture is that every unit interval graph is $e$-positive.
\citet{Hik24X} proved this equivalent form.
For Schur positivity,
\citet{Sta98} posed the conjecture that all claw-free graphs are Schur positive,
and attributed it to Gasharov.
It was disproved independently by \citet{Pra26X} and \citet{MM26X},
who found the same pair of counterexamples.

For individual graph families that are not unit interval,
\citet{LY21} proved~$e$-positivity for generalized pyramids,
\citet{TV26-HatChain} for hat chains,
and \citet{CHW26} for clocks.
These graphs are therefore Schur positive as well.
Of particular interest here are graphs that are Schur positive but not~$e$-positive.
Such graphs remain much less understood.
One such infinite family consists of the three-pendant generalized nets~$\mathrm{GN}_{n,3}$:
\citet{FKKMMT21} proved that they are not~$e$-positive,
and \citet{Sv25} later established Schur positivity for the larger family~$\mathrm{GN}_{n,m}$.
In contrast to these positivity results,
\citet{WW23-DAM} exhibited families of wheels and windmills that are not Schur positive,
and \citet{LLYZ25} proved the same for a family of squid graphs.

Necessary conditions for~$e$-positivity include that every $e$-positive graph contains a connected partition of every type, see \citet{Wol97D};
and that for Schur positivity include the \emph{niceness}:
if~$G$ contains a stable partition of type~$\lambda$,
then it contains a stable partition of any type that is dominated by~$\lambda$, 
see \citet{Sta98}.
\citet{LLYZ25} introduced the \emph{strongly nice property},
which is stronger than niceness and weaker than Schur positivity.
\citet{WW20} gave a complementary tool: a formula for Schur coefficients in terms of special ribbon tabloids.
Building on this formula, \citet{Sv25} introduced \emph{special rim hook $G$-tabloids}
and constructed sign-reversing maps on them to derive recurrences for Schur coefficients of generalized nets.

We study the distinction between the two positivity notions for spider graphs.
For a partition $\lambda=\lambda_1\lambda_2\dotsm\vdash n-1$,
the \emph{spider}~$S(\lambda)$ is the tree on~$n$ vertices obtained by identifying one endpoint of each of the paths~$P_{\lambda_i+1}$.
We also write~$S(a_1,\dots,a_r)$ for any sequence of positive leg lengths;
permuting the leg lengths does not change the isomorphism type.
See \cref{fig:Sa21.Sa41} for the spiders~$S(a,2,1)$ and~$S(a,4,1)$.

\begin{figure}[ht]
\begin{tikzpicture}[scale=.3]
\coordinate (0) at (0, 0);
\coordinate (1) at (2,0);
\coordinate (2) at (4,0); 
\coordinate (a1) at (4.5, 0); 
\coordinate (a2) at (5,0); 
\coordinate (a3) at (5.5,0); 
\coordinate (3) at (6,0); 
\coordinate (4) at (8,0); 
\coordinate (5) at (10,0); 
\coordinate (x) at (10,2); 
\coordinate (6) at (12,0);
\coordinate (7) at (14,0);
\draw[edge] (0)--(1)--(2);
\draw[edge] (3)--(4)--(5)--(6)--(7);
\draw[edge] (x)--(5);
\foreach \e in {0,1,4,5,x,6,7}   
\fill (\e) circle(.25);
\foreach \e in {a1, a2, a3}
\fill (\e) circle(.06);
\draw(a2) node[below=5pt]{length~$a$};
\draw(6) node[below=5pt]{length~$2$};

\begin{scope}[xshift=20cm]
\coordinate (0) at (0, 0);
\coordinate (1) at (2,0);
\coordinate (2) at (4,0); 
\coordinate (a1) at (4.5, 0); 
\coordinate (a2) at (5,0); 
\coordinate (a3) at (5.5,0); 
\coordinate (3) at (6,0); 
\coordinate (4) at (8,0); 
\coordinate (5) at (10,0); 
\coordinate (x) at (10,2); 
\coordinate (6) at (12,0);
\coordinate (7) at (14,0);
\coordinate (8) at (16,0);
\coordinate (9) at (18,0);
\draw[edge] (0)--(1)--(2);
\draw[edge] (3)--(4)--(5)--(6)--(7)--(8)--(9);
\draw[edge] (x)--(5);
\foreach \e in {0,1,4,5,x,6,7,8,9}   
\fill (\e) circle(.25);
\foreach \e in {a1, a2, a3}
\fill (\e) circle(.06);
\draw(a2) node[below=5pt]{length~$a$};
\draw(7) node[below=5pt]{length~$4$};
\end{scope}
\end{tikzpicture}
\caption{The spiders~$S(a,2,1)$ and~$S(a,4,1)$.}\label{fig:Sa21.Sa41}
\end{figure}

Spiders also connect this problem to the study of~$e$-positive trees.
\citet[Lemma~12]{DSv20} proved that if a tree~$T$ is~$e$-positive
and~$v$ is a vertex of the maximum degree~$\Delta(T)$,
then the spider whose leg lengths are the component orders of~$T-v$ has connected partitions of all types.
They further conjectured that no tree with~$\Delta\ge 4$ is~$e$-positive.
This has been confirmed for~$\Delta\ge 6$ by \citet{Zhe22},
and for~$\Delta=5$ and for spiders with~$\Delta=4$ by \citet{Tom26}.
The classification of~$e$-positive or Schur-positive trees with~$\Delta=3$ remains largely open.

Our choice of study objects is also motivated by structural simplicity.
Among trees that are not paths, three-legged spiders have the simplest branching structure:
a single branching vertex of degree~$3$.
Taking one leg to have the minimum length~$1$ gives the class~$S(a,b,1)$.
Fixing a second leg at successively larger lengths then leads to the one-parameter families
$S(a,1,1)$,~$S(a,2,1)$,~$\dots$.
For the one-parameter spider families~$S(a,1,1)$ and~$S(a,3,1)$,
\citet{WW23-DAM} showed that
\begin{itemize}
\item
$S(a,1,1)$ is~$e$-positive if and only if~$a=2$,
and is Schur positive but not~$e$-positive if and only if $a\in\{4,6,8,10,12\}$.
\item
$S(a,3,1)$ is~$e$-positive if and only if~$a\in\{2, 4\}$,
and is Schur positive but not~$e$-positive if and only if $a\in\{6,8,\dots,22\}$.
\end{itemize}
The first unresolved cases in this sequence were~$S(a,2,1)$ and~$S(a,4,1)$:
their Schur positivity had not been settled before this work,
although their~$e$-positivity classification was already known.
\citet{WW23-DAM} showed that
\begin{itemize}
\item
$S(a,2,1)$ is~$e$-positive if and only if~$a\in\{1,3,6\}$.
\item
$S(a,4,1)$ is~$e$-positive if and only if~$a\in\{3,5,8,10,12,13,15,20\}$.
\end{itemize}
In this paper, 
we prove that~$S(a,2,1)$ and~$S(a,4,1)$ are Schur positive for every integer~$a\ge1$;
see \cref{thm:spos:Sa21,thm:spos:Sa41}.
This completes the~$e$- and Schur-positivity classification of both families
and gives two infinite families of Schur-positive graphs that are not~$e$-positive.

Our proof lifts the composition-indexed path formula of \citet[Theorem~C.4]{SW16} at~$t=1$
to the algebra~$\mathrm{NSym}$ of noncommutative symmetric functions.
Its powers-of-two ribbon expansion (\cref{thm:path.ribbon}) and the Littlewood--Richardson rule
give a common Schur-coefficient formula for the spider reduction (\cref{lem:s-coeff:E}).
Weighted word maps in \cref{lem:multi-injection:j,lem:k=5} control the negative terms;
\cref{lem:finite-prefix} reduces the remaining comparisons to degree~$10$.

This~$\mathrm{NSym}$ realization is distinct from the chromatic symmetric function in noncommuting variables introduced by \citet{GS01},
which takes values in the algebra~$\mathrm{NCSym}$ of symmetric functions in noncommuting variables.
\citet{Cam24X} constructed a graph-wide lift taking values in~%
$\mathrm{NSym}$;
the problem-specific spider lifts used here are independent of that construction.

\citet[Section~1]{WZ25} described an earlier version of this paper as
``an early try of the composition method towards Schur positivity.''
Their method retains~$\mathrm{NSym}$ as theoretical support for composition-indexed~$e_I$-expansions;
later applications by \citet{CHW26,QTW26,Wang25} work directly with these expansions.
We retain~$\mathrm{NSym}$ because its ribbon basis is essential to our Schur-positivity arguments.
\citet[Theorem~9.9]{WW26SchurSpiders} later used the Schur positivity of~$S(a,2,1)$
as the base case of an induction proving that every spider~$S(a,b,2)$ is Schur positive.

After recalling the required background,
we develop a common path-based reduction for spiders~$S(a,b,1)$ in \cref{sec:NSym.Xpath}
and apply it to the two families in \cref{sec:proof.s+}.
\section{Preliminaries}

We use~$\mathbb Z$,~$\mathbb Q$, and~$\mathbb R$ for the integers, rationals, and reals,
with subscripts indicating restrictions, as in~$\mathbb Z_{\ge1}$.

\subsection{Words, compositions, partitions, and ribbons}

We use standard terminology for words,
factors,
prefixes,
suffixes,
and concatenation as in \citet{Lot05B}.
For a word $w=w_1w_2\dotsm w_l$,
write~$\ell(w)=l$ and~%
$\last(w)=w_l$.
For a word~$w$ over~$\mathbb Z_{\ge 1}$,
let~$m_i(w)$ denote the multiplicity of the letter~$i$ in~$w$,
and call the sequence $m_1(w)m_2(w)\dotsm$ its \emph{content}.
Write $\overline{w}=w_lw_{l-1}\dotsm w_1$ for the reversal of~$w$,
and use juxtaposition for concatenation.

A \emph{composition} of~$n$ is a word over the set of positive integers whose letters sum to~$n$,
denoted
$I
=i_1i_2\dotsm i_l
\vDash n$.
We regard~$0$ as having a unique composition,
denoted~$\emptyset$.
The \emph{sign} of~$I$ is~$(-1)^{n-l}$,
denoted~$\varepsilon^I$.
A \emph{refinement} of~$I$ replaces each part by a composition of that part.
We write~$I\preceq J$ when~$J$ refines~$I$, calling~$\preceq$ the \emph{reverse refinement order}.
For any nonempty compositions~$I=i_1\dotsm i_s$ and~$J=j_1\dotsm j_t$,
their \emph{concatenation}~$I\!J$ 
and \emph{near concatenation}~$I\triangleright J$ are the compositions
\[
I\!J
=(i_1,\,
\dots,\,
i_s,\,
j_1,\,
\dots,\,
j_t)
\quad\text{and}\quad
I\triangleright J
=(i_1,\,
\dots,\,
i_{s-1},\,
i_s+j_1,\,
j_2,\,
\dots,\,
j_t).
\]

A \emph{partition}\index{partition} of~$n$ is a multiset of positive integers with sum~$n$,
written~$\lambda=\lambda_1\lambda_2\dotsm\vdash n$ with~$\lambda_1\ge\lambda_2\ge\dots\ge1$.
We use French notation for Young diagrams.
A \emph{ribbon} is a connected skew shape with no~$2\times2$ block of boxes.
Throughout this paper, we identify a composition with the ribbon
whose row lengths, from top to bottom, are its parts.
We write~$\mathrm{sh}(I)$ for the skew shape corresponding to a composition~$I$.
For example,
\[
\mathrm{sh}(12)
=
\begin{tikzpicture}[baseline={([yshift=-.25ex]current bounding box.center)},scale=.28]
\draw (0,0)--(2,0)--(2,1)--(1,1)--(1,2)--(0,2)--cycle;
\draw (1,0)--(1,1);
\draw (0,1)--(1,1);
\end{tikzpicture}
=21
\qquad\text{and}\qquad
\mathrm{sh}(21)
=
\begin{tikzpicture}[baseline={([yshift=-.25ex]current bounding box.center)},scale=.28]
\draw (0,1)--(0,2)--(2,2)--(2,0)--(1,0)--(1,1)--cycle;
\draw (1,1)--(1,2);
\draw (1,1)--(2,1);
\end{tikzpicture}
=22/1.
\]
Following \citet{Mac1915B},
the \emph{conjugate} of a composition~$I$ is the composition formed by the column heights of its ribbon \emph{from right to left},
denoted~$I^\sim$.
Reversal commutes with conjugation, so~$\overline{I}^\sim=\overline{I^\sim}$
reads the column heights from left to right.

A \emph{run} is a maximal strictly increasing factor of a word.
Every word~$w$ has a unique run factorization~%
$w=\xi_1\dotsm\xi_r$.
We define its \emph{run type}~$\tau_w$ to be the composition of run lengths,
namely,
$\tau_w=\ell(\xi_1)\dotsm\ell(\xi_r)$.
For a word~$w$ over~$\mathbb Z_{\ge 1}$,
its \emph{shape} is the ribbon shape~$\mathrm{sh}(\tau_w^\sim)$.

\subsection{Commutative and noncommutative symmetric functions}\label{sec:Sym}\label{sec:NSym}
For background on commutative symmetric functions,
see \citet{Mac95B} and \citet[Chapter~7]{Sta23B}.
Let~$\mathrm{Sym}$ be the algebra of symmetric functions over~$\mathbb Q$.
We use the standard notation~$e_\lambda$,~$p_\lambda$, and~$s_\lambda$
for the elementary, power-sum, and Schur basis elements of degree~$n$.
For a skew shape~$\lambda/\mu$,
we write~$s_{\lambda/\mu}$ for the corresponding skew Schur function.
We call a symmetric function~\emph{$e$-/Schur positive}
if all coefficients in its~$e$-/Schur expansion
are nonnegative.
A \emph{Yamanouchi word} is a word~$y$ such that
\[
m_i(y')\ge m_{i+1}(y')
\quad\text{for every prefix $y'$ of $y$
and every positive integer $i$.}
\]
These are the \emph{Yamanouchi inequalities}, also called the \emph{ballot inequalities}.
The content of a Yamanouchi word is therefore a partition.

\begin{theorem}[Littlewood and Richardson]\label{thm:Littlewood--Richardson}
Let~$I\vDash n$ and~$\nu\vdash n$.
Then the coefficient~$[s_\nu]s_{\mathrm{sh}(I)}$
equals the number of Yamanouchi words of shape~$\mathrm{sh}(I)$ and content~$\nu$.
\end{theorem}

For a filling of a skew Young diagram in French notation,
its \emph{reading word} is the word obtained by reading the entries row by row
from bottom to top and from right to left within each row.
This is the reading convention used in the Littlewood--Richardson rule.

For an introduction to noncommutative symmetric functions and their basic properties,
see \citet{GKLLRT95}.
The \emph{algebra of noncommutative symmetric functions} is the free associative algebra
\[
\mathrm{NSym}
=\mathbb Q\langle\Lambda_1,\Lambda_2,\dots\rangle
\]
generated over~$\mathbb Q$ by an infinite sequence~$\{\Lambda_k\}_{k\ge 1}$ of indeterminates,
where~$\Lambda_0=1$.
The generators~$\Lambda_n$ are the noncommutative elementary symmetric functions.
Set
$\lambda(t)
=\sum_{n\ge0}\Lambda_n t^n$,
where~$t$ commutes with every~$\Lambda_n$.
The \emph{power sum symmetric functions of the first kind},
denoted~$\Psi_n$,
are defined by
\[
\sum_{n\ge 1}\Psi_n t^{n-1}
=\lambda(-t)
\frac{d \lambda(-t)^{-1}}{dt}.
\]

For any composition~$I=i_1 i_2\dotsm$,
define
$\Lambda^I
=\Lambda_{i_1}
\Lambda_{i_2}\dotsm$,
and 
$\Psi^I
=\Psi_{i_1}
\Psi_{i_2}\dotsm$.
Let~$R_I$ denote the \emph{noncommutative ribbon Schur function} indexed by~$I$.
There is a unique graded algebra homomorphism
\[
\rho\colon\mathrm{NSym}\to\mathrm{Sym}
\]
such that~$\rho(\Lambda_n)=e_n$ for every~$n\ge1$.
We call~$\rho$ the \emph{abelianization map}.
It follows that
\[
\rho(\Psi_n)=p_n
\quad\text{and}\quad
\rho(R_I)=s_{\mathrm{sh}(I)}.
\]
For a composition~$I=i_1\dotsm i_s$, define
$e_I=e_{i_1}\dotsm e_{i_s}$
and
$p_I=p_{i_1}\dotsm p_{i_s}$.
Then
\[
\rho(\Lambda^I)=e_I
\quad\text{and}\quad
\rho(\Psi^I)=p_I.
\]
When~$\rho(F)=f$ for some~$F\in\mathrm{NSym}$ and~$f\in\mathrm{Sym}$,
we call~$F$ a \emph{noncommutative analog} of~$f$.
By \cref{thm:Littlewood--Richardson}, the coefficient~$[s_\nu]\rho(R_I)$
is the number of Yamanouchi words of shape~%
$\mathrm{sh}(I)$ and content~$\nu$.
Ribbon Schur functions form a linear basis of~$\mathrm{NSym}$
and their multiplication satisfies a simple rule,
see \citet[Proposition 3.13]{GKLLRT95}.

\begin{proposition}[\citeauthor{GKLLRT95}]\label{prop:RI.RJ}
We have $R_I R_J
=R_{I\!J}+R_{I\triangleright J}$
for any nonempty compositions~$I$ and~$J$.
For the empty composition,~$R_{\emptyset}=1$.
\end{proposition}

We will also use some transition formulas between bases,
see \citet[Proposition~4.15 and Note~4.21]{GKLLRT95}.
To state them,
we need some additional notation.
Let $I=i_1\dotsm i_s$ and $J=j_1\dotsm j_t$ be compositions of~$n$.
Assume that~$J\succeq I$.
More precisely,
suppose that
\[
i_1=j_{k_0+1}+\dots+j_{k_1},\quad
\dots,\quad
i_s=j_{k_{s-1}+1}+\dots+j_{k_s},
\]
for some indices $0=k_0<\dots<k_s=t$.
We define
\[
\fp(J,I)=\prod_{r=1}^s j_{k_{r-1}+1}
\]
to be the product of the \emph{first parts of the blocks of~$J$ relative to~$I$}.
We now state the basis transition formulas used below.

\begin{proposition}[\citeauthor{GKLLRT95}]\label{prop:transition:NSym}
For any composition~$I=i_1\dotsm i_s$, we have
\[
\Lambda^I
=\sum_{J\succeq \overline{I}^\sim}R_J,
\quad
\Psi^I
=\sum_{J\succeq I}
\fp(J,I)
\varepsilon^J \Lambda^J, 
\quad\text{and}\quad
R_I
=\sum_{J\preceq\overline{I}^\sim}
(-1)^{\ell(\overline{I}^\sim)-\ell(J)}\Lambda^J.
\]
\end{proposition}

\section[Path expansions and a Schur-coefficient reduction]{Path expansions and the Schur-coefficient reduction for spiders}\label{sec:NSym.Xpath}

This section develops the common reduction used in the two Schur-positivity proofs in \cref{sec:proof.s+}.
We first separate an~$e$-positive summand from the chromatic symmetric function of~$S(a,b,1)$.
We then use ribbon expansions in~$\mathrm{NSym}$ and the Littlewood--Richardson rule
to express the Schur coefficients of the remaining summand in terms of weighted Yamanouchi words,
obtaining the formula in \cref{lem:s-coeff:E}.

\subsection{The path lift and the positive summand}\label{sec:path-lift}

Let~$a,b\ge1$ and set~$n=a+b+2$.
\citet{Zhe22} obtained the following reduction formula for spiders~$S(a,b,1)$:
\begin{equation}\label{X.Sab1}
X_{S(a,b,1)}
=X_{P_{a+b+2}}
+e_1 X_{P_{n-1}}
-X_{P_{a+1}}X_{P_{b+1}}.
\end{equation}
For any~$n\ge 2$, there is a unique pair~$(A_{n-1},\, B_n)$ of symmetric functions such that 
\begin{equation}\creflabel[def]{def:path.AB}
X_{P_n}=e_1 A_{n-1}+ B_n,
\end{equation}
with~$B_n\in\mathbb Q[e_2,e_3,\dots]$.
Uniqueness follows from the algebraic independence of the elementary generators.
Define
\[
D_{n,k}
=X_{P_n}-X_{P_{n-k}}A_k
\quad\text{and}\quad
E_{n,k}
=X_{P_n}-X_{P_{n-k}}B_k.
\]
Then by \cref{X.Sab1},
\begin{equation}\label{X.Sab1:DE}
X_{S(a,b,1)}
=X_{P_n}
+e_1 X_{P_{n-1}}
-X_{P_{a+1}}(e_1 A_{b}+B_{b+1})
=e_1 D_{n-1,\,b}+E_{n,\,b+1}.
\end{equation}
We will deduce the~$e$-positivity of~$D_{n,k}$ in \cref{cor:epos.A}, 
and prove the Schur positivity of~$E_{n,3}$ and~$E_{n,5}$ in
\cref{thm:spos:Sa21,thm:spos:Sa41}, respectively.

For a composition~$I=i_1\dotsm i_s$,
define
\[
w_I=i_1(i_2-1)(i_3-1)\dotsm(i_s-1).
\]
\citet[Theorem~C.4]{SW16} obtained the following formula at~$t=1$ using Smirnov words;
we give a proof in~$\mathrm{NSym}$.

\begin{theorem}[\citeauthor{SW16}]\label{thm:SW.path}
For any~$n\ge1$, $X_{P_n}
=\sum_{I\vDash n}w_I e_I$.
Equivalently, the chromatic symmetric function~$X_{P_n}$ has a noncommutative analog
$\widetilde{X}_{P_n}
=\sum_{I\vDash n}w_I\Lambda^I$.
\end{theorem}

\begin{proof}
Since~$\rho(\Lambda^I)=e_I$, it suffices to prove the noncommutative formula.
\citet[Theorem~2.5]{Sta95} showed that
\[
X_G
=\sum_{E'\subseteq E(G)}(-1)^{\abs{E'}}p_{\tau(E')},
\]
where~$\tau(E')$ is the partition consisting of
the component orders of the spanning subgraph~$(V(G),E')$.
For~$G=P_n$,
reading off the component orders gives a bijection from the power set of~$E(P_n)$ to the set of compositions of~$n$.
Therefore,
$X_{P_n}
=\sum_{I\vDash n}\varepsilon^I p_I$.
Since~$\rho(\Psi^I)=p_I$,
the formula above gives the noncommutative analog 
$\widetilde{X}_{P_n}
=\sum_{I\vDash n}\varepsilon^I \Psi^I$.
By \cref{prop:transition:NSym},
we obtain
\[
\widetilde{X}_{P_n}
=\sum_{I\vDash n}
\varepsilon^I
\sum_{J\succeq I}
\fp(J,I)
\varepsilon^J \Lambda^J
=\sum_{J\vDash n}
\brk4{
\sum_{I\preceq J}
\varepsilon^{I\!J}
\fp(J,I)}
\Lambda^J.
\]
Let $J=j_1\dotsm j_t\vDash n$.
For any $I=i_1\dotsm i_s\preceq J$, 
let~$j_{k_h}$ be the first part of the~$h$th block of~$J$ relative to~$I$.
Then $1=k_1<\dots<k_s\le t$,
$\varepsilon^{I\!J}
=(-1)^{t-s}$,
and $\fp(J,I)
=j_{k_1}\dotsm j_{k_s}$.
Therefore,
\[
\brk[s]1{\Lambda^J}\widetilde{X}_{P_n}
=\sum_{I\preceq J}
\varepsilon^{I\!J}
\fp(J,I)
=j_1\sum_{s=1}^t
\sum_{2\le k_2<\dots<k_s\le t}
(-1)^{t-s}j_{k_2}\dots j_{k_s}
=w_J.\qedhere
\]
\end{proof}

We now apply the path expansion to~$D_{n,k}$ in the spider reduction~\cref{X.Sab1:DE}.

\begin{corollary}\label{cor:epos.A}
For any $1\le k\le n-1$,
the symmetric function~$D_{n,k}$ admits a noncommutative analog
whose expansion in the~$\Lambda^I$ basis has nonnegative coefficients.
Consequently,~$D_{n,k}$ is~$e$-positive.
\end{corollary}

\begin{proof}
In \cref{thm:SW.path}, a part~$i_j=1$ with~$j\ge2$ forces~$w_I=0$.
Hence, among the terms with nonzero coefficient, a factor~$e_1$ occurs in the commutative image
exactly when~$i_1=1$.
Therefore, 
the symmetric function~$e_1 A_{n-1}$
has a noncommutative analog~$\Lambda^1\widetilde{A}_{n-1}$ such that 
\[
\widetilde{A}_{n-1}
=
\sum_{1i_2\dotsm i_s\vDash n}
(i_2-1)\dotsm(i_s-1)\Lambda^{i_2\dotsm i_s}.
\]
By definition, the symmetric function~$D_{n,k}$ has the noncommutative analog
\begin{align*}
\widetilde{D}_{n,k}
&=\widetilde{X}_{P_n}
-\widetilde{X}_{P_{n-k}}\widetilde{A}_k
=
\sum_{I\vDash n}
w_I\Lambda^I
-
\sum_{I\vDash n-k}
w_I\Lambda^I
\sum_{J=j_1\dotsm j_t\vDash k}
(j_1-1)\dotsm(j_t-1)\Lambda^J
\\
&=
\sum_{I=i_1\dotsm i_s\vDash n, \
i_1+\dotsm+i_j\ne n-k\text{ for }1\le j\le s}
w_I\Lambda^I.
\end{align*}
Hence the expansion of~$\widetilde{D}_{n,k}$ in the~$\Lambda^I$ basis
has nonnegative coefficients.
Its abelianization~$D_{n,k}$ is therefore~$e$-positive.
\end{proof}

\subsection{The ribbon expansion of the path lift}\label{sec:path-ribbon}

For~$n\ge 0$, let~$\Cge{n}$ and~$\Cle{n}$ be the sets of
compositions of~$n$ all of whose parts are at least~$2$ and at most~$2$,
respectively. 
For example,
\[
\Cge{4}=\{4,\,22\}
\quad\text{and}\quad
\Cle{4}=\{22,\,211,\,121,\,112,\,1^4\}.
\]
We use the convention $\Cge{0}=\Cle{0}=\{\emptyset\}$.
For a set~$\mathcal A$ of compositions and compositions~$K,L$, write
\[
K\mathcal A=\{K\!I\colon I\in\mathcal A\},
\qquad
\mathcal A L=\{I\!L\colon I\in\mathcal A\},
\qquad\text{and}\quad
K\mathcal A L=\{K\!I\!L\colon I\in\mathcal A\}.
\]

\begin{theorem}\label{thm:path.ribbon}
For any~$n\ge1$, the analog~$\widetilde{X}_{P_n}$ defined in \cref{thm:SW.path}
satisfies
\[
\widetilde{X}_{P_n}
=\sum_{I\in\Cle{n-1}}
2^{m_1(I)}R_{I1}
=
\sum_{J\in \mathcal R_n}
R_J,
\]
where~$m_1(I)$ is the number of parts~$1$ in~$I$,
and~$\mathcal R_n$ is the set of ribbons with~$n$ boxes labeled by~$0$ or~$1$
such that every horizontal adjacency, read from left to right, has labels~$10$
and the last row is a single box labeled by~$1$.
\end{theorem}

\begin{proof}
For a binary word~$b=b_1\dotsm b_n$ ending in~$1$,
cut between~$b_k$ and~$b_{k+1}$ whenever~$b_kb_{k+1}=10$,
and let~$I(b)$ be the composition of the block lengths.
Each block is weakly increasing and ends in~$1$,
and every block except the first begins with~$0$.
For fixed~$I=i_1\dotsm i_s\vDash n$,
there are~$i_1$ choices for the first block and~$i_j-1$ choices for each subsequent block~$j$.
Thus exactly~$w_I$ words satisfy~$I(b)=I$,
and \cref{thm:SW.path} gives
\[
\widetilde{X}_{P_n}
=\sum_{b\in\{0,1\}^n,\ b_n=1}\Lambda^{I(b)}.
\]

Since~$\Lambda_r=R_{1^r}$,
\cref{prop:RI.RJ} expands~$\Lambda^{I(b)}$ by joining the corresponding vertical ribbons
vertically or horizontally at each block boundary.
For a fixed~$b$, each resulting ribbon occurs once,
since its adjacencies determine the choice at every boundary.
For a fixed ribbon~$J$, label its boxes by~$b_1,\dots,b_n$,
reading rows from top to bottom and each row from left to right.
The coefficient~$[R_J]\widetilde{X}_{P_n}$ therefore counts exactly the words~$b$ ending in~$1$
for which every horizontal adjacency is labeled~$10$;
vertical adjacencies impose no restriction.

Two successive horizontal adjacencies would force their common letter to be both~$0$ and~$1$,
so every row has length at most~$2$.
The last row must have length~$1$, since the word ends in~$1$.
Thus these labeled ribbons are precisely the members of~$\mathcal R_n$,
and their sum is~$\widetilde{X}_{P_n}$.
For each~$I\in\Cle{n-1}$, a ribbon of shape~$I1$ has fixed labels~$10$ on its two-box rows
and~$1$ on its final box,
while each of its other~$m_1(I)$ singleton rows can independently be labeled~$0$ or~$1$.
Thus exactly~$2^{m_1(I)}$ members of~$\mathcal R_n$ have shape~$I1$.
Grouping the labeled ribbons by shape proves the formula.
\end{proof}

For an alternative proof by a sign-reversing involution, see \cref{sec:path-involution}.

\subsection{The complementary path term}\label{sec:path-complement}

To construct a noncommutative analog of the symmetric function~$B_n$ introduced in \cref{def:path.AB},
we will need \cref{lem:B,lem:m1:ell.sim}.

\begin{lemma}\label{lem:B}
Let~$n\ge 2$. Then 
\[
\sum_{I=i_1\dotsm i_s\in\Cge{n}}
(i_1-1)\dotsm(i_s-1)
=2^{n-2}
\quad\text{and}\quad
\sum_{I\in\Cge{n}}w_I
=\begin{cases}
2,&\text{if $n=2$},\\[3pt]
3\cdotp 2^{n-3},&\text{if $n\ge 3$}.
\end{cases}
\]
\end{lemma}
\begin{proof}
For~$n\ge 2$,
define
\begin{equation}\label{pf:def:fg}
f_n
=\sum_{I=i_1\dotsm i_s\in\Cge{n}}
(i_1-1)\dotsm(i_s-1)
\quad\text{and}\quad
g_n
=\sum_{I\in\Cge{n}}w_I.
\end{equation}
Set~$f_0=1$ and~$f_1=0$, corresponding to the empty composition and the absence of compositions in~$\Cge{1}$.
Then~$f_2=1$ and~$g_2=2$. Let~$n\ge 3$.
By \cref{pf:def:fg},
\[
f_n
=\sum_{c_1=2}^{n}
(c_1-1)
f_{n-c_1}.
\]
Writing~$F(z)=\sum_{n\ge0}f_nz^n$, the recurrence gives
\[
F(z)=1+F(z)\sum_{r\ge2}(r-1)z^r
=1+\frac{z^2}{(1-z)^2}F(z).
\]
Consequently,
\[
\sum_{n\ge2}f_nz^n=F(z)-1
=\frac{z^2}{1-2z}
=\sum_{n\ge2}2^{n-2}z^n.
\]
This proves the first desired identity~%
$f_n=2^{n-2}$.

With the convention that an empty sum equals~$0$,
separating the one-part compositions in the definition of~$g_n$
and writing $i_1=1+(i_1-1)$ for each composition with at least two parts,
we obtain
\begin{multline*}
g_n
=n+\sum_{I\in\Cge{n},\ \ell(I)\ge 2}w_I
=n
+\sum_{I=i_1\dotsm i_s\in\Cge{n},\ s\ge 2}
\prod_{r=2}^s(i_r-1)
+\sum_{I=i_1\dotsm i_s\in\Cge{n},\ s\ge 2}
\prod_{r=1}^s(i_r-1)\\
=n
+\sum_{i_1=2}^{n-2}f_{n-i_1}
+(f_n-n+1)
=
\sum_{j=0}^{n-4}2^j+2^{n-2}+1
=3\cdotp 2^{n-3}.\qedhere
\end{multline*}
\end{proof}

To prove the following path decomposition,
we express the number~$m_1(I)$ of parts equal to~$1$
in terms of the conjugate composition~$I^\sim$.

\begin{lemma}\label{lem:m1:ell.sim}
Let~$n\ge1$, and let~$I\vDash n$ be a composition all of whose parts are~$1$ or~$2$. Then 
\[
m_1(I)
=n+2-2\ell(I^\sim).
\]
\end{lemma}

\begin{proof}
Each part of~$I$ equal to~$2$ starts a new column.
Then the desired formula follows from the identities
$\ell(I^\sim)=m_2(I)+1$ and
$n=m_1(I)+2m_2(I)$.
\end{proof}

For example,
the Yamanouchi word~$y=121213$ has run factorization
$y=(12)(12)(13)$,
so~$\tau_y=222$.
The composition $I=\tau_y^\sim=1221$ and its conjugate $I^\sim=\tau_y=222$ have shapes
\[
\mathrm{sh}(I)
=
\begin{tikzpicture}[baseline={([yshift=-.25ex]current bounding box.center)},scale=.28]
\foreach \x/\y in {0/3,0/2,1/2,1/1,2/1,2/0}
\draw (\x,\y) rectangle +(1,1);
\end{tikzpicture}
\qquad\text{and}\qquad
\mathrm{sh}(I^\sim)
=
\begin{tikzpicture}[baseline={([yshift=-.25ex]current bounding box.center)},scale=.28]
\foreach \x/\y in {0/2,1/2,1/1,2/1,2/0,3/0}
\draw (\x,\y) rectangle +(1,1);
\end{tikzpicture}
\]
respectively.
Then \cref{lem:m1:ell.sim} gives
$m_1(I)=6+2-2\ell(222)=2$.

\begin{lemma}\label{lem:Y:B.path}
We have~$B_2=\rho(2R_{11})$,
and for~$n\ge 3$ we have~$B_n=\rho(\widetilde{B}_n)$,
where
\[
\widetilde{B}_n
=
\sum_{I=2\dotsm\in\Cle{n-2}}
2^{m_1(I)+1}R_{1I1}
+3
\sum_{I=1\dotsm\in\Cle{n-2}}
2^{m_1(I)-1}R_{1I1}.
\]
Here the notation~$x\dotsm$ stands for a composition whose first part is~$x$.
\end{lemma}
\begin{proof}
Since~$\overline{11}^\sim=2$ and~$2$ has no proper coarsening,
\cref{prop:transition:NSym} gives~$\rho(R_{11})=e_2$.
Therefore,
\[
\rho(2R_{11})
=2e_2
=B_2.
\]
Now let~$n\ge 3$.
In view of \cref{thm:SW.path},
the compositions with~$i_1=1$ contribute to~$e_1 A_{n-1}$,
whereas the compositions with~$i_1\ge2$ contribute to~$B_n$.
Therefore,~%
$B_n$ admits the noncommutative analog
\begin{align*}
\widetilde{B}_n
&=
\sum_{I=i_1\dotsm i_s\vDash n, \
i_1\ge 2}
w_I\Lambda^I
=
\sum_{I\in\Cge{n}}
\sum_{J\succeq \overline{I}^\sim}
w_I
R_J, 
\end{align*}
where for deriving the last equality we used \cref{prop:transition:NSym}.
We now exchange the summation order for this double sum.
Let $I=i_1\dotsm i_s\in\Cge{n}$.
Reading the column heights of~$\overline{I}$ from right to left,
we find that~$\overline I^\sim$ has every part at most~$2$
and starts and ends with a part~$1$.
In other words,~$\overline{I}^\sim\in 1\Cle{n-2}1$.
Since refining a part~$2$ can only replace it by~$11$
and a part~$1$ cannot be refined,
the set~$1\Cle{n-2}1$ forms an order filter under refinement~$\preceq$.
Thus every composition~$J\succeq\overline I^\sim$ also belongs to~%
$1\Cle{n-2}1$.
Since conjugation is an involution that reverses refinement, 
\[
J\succeq \overline{I}^\sim
\iff
\overline I=(\overline I^\sim)^\sim\succeq J^\sim.
\]
Since reversal is an involution that preserves refinement,
\[
J\succeq \overline{I}^\sim
\iff
I\succeq \overline{J^\sim}.
\]
Therefore, by exchanging the summation order, we obtain
\begin{align*}
\widetilde{B}_n
&=\sum_{J\in 1\Cle{n-2}1, \
I\in\Cge{n}, \ I\succeq \overline{J^\sim}}
w_I
R_J.
\end{align*}
Fix~$J\in 1\Cle{n-2}1$. 
Then the condition for~$I$ says~$I$ refines the column heights~$t_i$ of~$J$.
It follows that
\[
[R_J]\widetilde{B}_n
=\sum_{\substack{
c_{1k}+c_{2k}+\dots=t_k,\,\forall\,k\\
c_{ij}\ge 2,\ \forall\, i,\,j}}
w_{(c_{11}, c_{21}, \dots,
c_{12}, c_{22}, \dots,
c_{1l}, c_{2l}, \dots)}
=g_{t_1}\prod_{i\ge 2}f_{t_i},
\]
where
\[
g_n
=\sum_{I\in\Cge{n}}w_I
\quad\text{and}\quad
f_n
=\sum_{I=i_1\dotsm i_s\in\Cge{n}}
(i_1-1)\dotsm(i_s-1).
\]
By \cref{lem:B,lem:m1:ell.sim},
we obtain
\[
[R_J]\widetilde{B}_n
=\begin{cases}
2 \prod_{i\ge 2}2^{t_i-2}
=2^{1+(n-t_1)-2(\ell(J^\sim)-1)}
=2^{m_1(J)-1},
&\text{if $t_1=2$},
\\[5pt]
3\cdotp 2^{t_1-3} \prod_{i\ge 2}2^{t_i-2}
=3\cdotp 2^{n-3-2(\ell(J^\sim)-1)}
=3\cdotp 2^{m_1(J)-3},
&\text{if $t_1\ge 3$}.
\end{cases}
\]
Writing~$J=1I1$ with~$I\in\Cle{n-2}$,
we have $m_1(J)=m_1(I)+2$.
Moreover,~$t_1=2$ if and only if~$I=2\dotsm$,
whereas~$t_1\ge3$ if and only if~$I=1\dotsm$.
Substituting these relations into the coefficient formula gives
\[
[R_{1I1}]\widetilde{B}_n
=\begin{cases}
2^{m_1(I)+1},&\text{if $I=2\dotsm$},\\[3pt]
3\cdotp 2^{m_1(I)-1},&\text{if $I=1\dotsm$}.
\end{cases}
\]
This proves the desired formula for~$\widetilde{B}_n$.
\end{proof}

\subsection{Weighted words and the common coefficient formula}\label{sec:weighted-words}

A \emph{multiset}~$A$ consists of a finite \emph{support}~$\underline A$
and a multiplicity function
$m_A\colon\underline A\to\mathbb R_{>0}$,
extended by zero outside~$\underline A$.
It has \emph{norm} 
\[
\norm{A}
=\sum_{a\in\underline A}m_A(a).
\]
We write~$a\in A$ when~$a\in\underline A$.
If $\underline A=\{a_1,\dots,a_x\}$ and~$m_i=m_A(a_i)$,
we use the compact notation
\[
A=\{m_1\cdotp a_1,\,
\dots,\,
m_x\cdotp a_x\}.
\]
More generally,
for any finite set~$S$ and any function~%
$g\colon S\to\mathbb R_{>0}$,
we write~$(S,g)$ for the multiset whose multiplicity function is~$g$ on~$S$
and~$0$ outside~$S$.
For~$r\in\mathbb R_{\ge0}$,
the \emph{scalar multiset}~$rA$ has multiplicity function~$r m_A$.
A multiset~$B$ is a \emph{sub-multiset} of~$A$ if
$m_B(u)\le m_A(u)$ for every~$u$.
The disjoint union and difference are defined pointwise by
\[
m_{A\sqcup B}(u)
=m_A(u)+m_B(u)
\quad\text{and}\quad
m_{A\backslash B}(u)
=\max\brk1{m_A(u)-m_B(u),\,0}.
\]

A \emph{multi-map}~$f$ from~$A$ to~$B$ is specified by a family~%
$(k_{ab})_{a\in\underline A,\,b\in\underline B}$
of nonnegative real numbers.
We regard~$f(a)$ as the multiset
\[
f(a)=\{k_{ab}\cdotp b\colon b\in\underline B,\ k_{ab}>0\}.
\]
The multi-map is \emph{mass-preserving} if
\[
\sum_{b\in\underline B}k_{ab}=1
\quad\text{for every }a\in\underline A.
\]
Its \emph{image}~$f(A)$ is the multiset determined by
\[
m_{f(A)}(b)
=\sum_{a\in\underline A}m_A(a)k_{ab}
\quad\text{for every }b\in\underline B.
\]
We say that~$f$ is a \emph{multi-injection} if
$m_{f(A)}(b)\le m_B(b)$ for every~$b\in\underline B$.

In order to show the Schur positivity of~$E_{n, k}$, 
we give a combinatorial interpretation for the ribbon coefficients of the noncommutative analog
\[
\widetilde{E}_{n,k}
=\widetilde{X}_{P_n}-\widetilde{X}_{P_{n-k}}\widetilde{B}_k, 
\]
in terms of multisets of Yamanouchi words.

For any partition~$\kappa\vdash n$,
let~$\mathcal{Y}(\kappa)$ be the set of Yamanouchi words~$y$ of content~$\kappa$
such that only the last run of~$y$ may have length~$1$.
Let~$\mathcal{Y}_\alpha(\kappa)$ denote the set of words~$y\in\mathcal{Y}(\kappa)$ whose run type~$\tau_y$ has prefix~$\alpha$.
For example,
\[
\mathcal{Y}(321)
=\{121213,\,121231,\,121312,\,123121\}
\quad\text{and}\quad
\mathcal{Y}_{22}(321)=\{121213,\,121312\}.
\]
Define the \emph{weight} of a Yamanouchi word~$y$ by
\[
\wt(y)
=2^{m_1(\tau_y^\sim)-1}.
\]
Let~$y\in\mathcal{Y}(\kappa)$ with length~$n$.
Then the last-run requirement for~$y$ is equivalent to~$\tau_y^\sim\in\Cle{n-1}1$.
By \cref{lem:m1:ell.sim},
\begin{equation}\label{wt.y}
\wt(y)
=2^{n+1-2\ell(\tau_y)}.
\end{equation}
Define the multiset
$M_\alpha(\kappa)=(\mathcal{Y}_\alpha(\kappa),\,\wt)$.

We allow the empty prefix, so that~$\mathcal Y_\emptyset(\kappa)=\mathcal Y(\kappa)$;
recall that~$\Cge{0}=\{\emptyset\}$.

\begin{proposition}\label{lem:s-coeff:E}
Let~$k\ge3$,~$n\ge 2k$, and~$\nu\vdash n$.
Then
\begin{align*}
[s_\nu]E_{n,k}
&=
\sum_{\alpha\in\Cge{k-1}}
\norm{M_{\alpha}(\nu)}
+\sum_{\alpha\in\Cge{k+1},\,\last(\alpha)\ge 3}
\norm{M_{\alpha}(\nu)}
+\frac{1}{2}
\sum_{\alpha\in\Cge{k-2},\,z\ge 4}
\norm{M_{\alpha z}(\nu)}\\
&\qquad
+\frac{5}{8}\sum_{i=3}^k
\sum_{\alpha\in\Cge{k-i},\,z\ge i+2}
\norm{M_{\alpha z}(\nu)}
-\smashoperator{
\sum_{\alpha\in\Cge{k},\,\last(\alpha)=2}}
\norm{M_{\alpha}(\nu)}
-\frac{1}{2}
\sum_{\alpha\in\Cge{k},\,\last(\alpha)\ge 3}
\norm{M_{\alpha}(\nu)}.
\end{align*}
\end{proposition}

\begin{proof}
Recall that $E_{n,k}=X_{P_n}-X_{P_{n-k}}B_k$.
By \cref{thm:path.ribbon,lem:Y:B.path},
\begin{align*}
\widetilde{X}_{P_n}
&=\sum_{I\in\Cle{n-1}}
2^{m_1(I)}R_{I1}, \quad\text{and}\\
\widetilde{X}_{P_{n-k}}\widetilde{B}_k
&=\sum_{\substack{
\alpha\in\Cle{n-k-1}\\
\beta=2\dotsm\in\Cle{k-2}}}
2^{m_1(\alpha)+m_1(\beta)+1}
R_{\alpha1}R_{1\beta1}
+3\sum_{\substack{
\alpha\in\Cle{n-k-1}\\
\beta=1\dotsm\in\Cle{k-2}}}
2^{m_1(\alpha)+m_1(\beta)-1}
R_{\alpha1}R_{1\beta1}, 
\end{align*}
in which $R_{\alpha1}R_{1\beta1}
=R_{\alpha11\beta1}+R_{\alpha2\beta1}$
by \cref{prop:RI.RJ}. Writing out the first letter of~$\beta$, we obtain
\begin{align*}
\widetilde{E}_{n,k}
&=
\sum_{I\in\Cle{n-1}}
2^{m_1(I)}R_{I1}
-\sum_{
\alpha\in\Cle{n-k-1}, \
\beta\in\Cle{k-4}}
2^{m_1(\alpha)+m_1(\beta)+1}
(R_{\alpha 22\beta 1}+R_{\alpha112\beta1})\\
&\quad{}
-3\sum_{
\alpha\in\Cle{n-k-1}, \
\beta\in\Cle{k-3}}
2^{m_1(\alpha)+m_1(\beta)}
(R_{\alpha21\beta1}+R_{\alpha111\beta1}).
\end{align*}
The first negative sum in this display is understood to be empty when~$k=3$.
Each of the four ribbon families in the negative sums also occurs in the positive sum.
The terms in the positive sum outside these four ribbon families form
\[
\widetilde{U}_{n,k}
=
\sum_{
\gamma\in\Cle{n-k-2}, \
\delta\in\Cle{k-1}}
2^{m_1(\gamma2\delta)}R_{\gamma2\delta1}
+\sum_{
\alpha\in\Cle{n-k-1}, \
\xi\in\Cle{k-3}}
2^{m_1(\alpha12\xi)}R_{\alpha12\xi1}.
\]

We first express each Schur coefficient of~$\rho(\widetilde{U}_{n,k})$
as a sum of weights in the multisets~$M_\theta(\nu)$,
then compute the cancellation among the four ribbon families,
and finally translate the resulting difference into multiset weights.

For any~$I\in\Cle{n-1}$ and~$\nu\vdash n$,
\cref{thm:Littlewood--Richardson} identifies
the number $[s_\nu]\rho(R_{I1})=[s_\nu]s_{\mathrm{sh}(I1)}$
with the Yamanouchi words~$y\in\mathcal Y(\nu)$ satisfying
$\tau_y=(I1)^\sim$.
For each such word,
$\wt(y)
=2^{m_1(\tau_y^\sim)-1}
=2^{m_1(I)}$.
Consequently, 
\[
2^{m_1(I)}[s_\nu]s_{\mathrm{sh}(I1)}
=
\sum_{
y\in\mathcal Y(\nu), \
\tau_y=I1^\sim}
\wt(y).
\]
It follows that 
\begin{align*}
[s_\nu]\rho(\widetilde{U}_{n,k})
&=
\sum_{
\gamma\in\Cle{n-k-2}, \
\delta\in\Cle{k-1}}
2^{m_1(\gamma2\delta)}
[s_\nu]
s_{\mathrm{sh}(\gamma2\delta1)}
+
\sum_{
\alpha\in\Cle{n-k-1}, \
\xi\in\Cle{k-3}}
2^{m_1(\alpha12\xi)}
[s_\nu]
s_{\mathrm{sh}(\alpha12\xi1)}
\\
&=
\sum_{
y\in\mathcal Y(\nu), \
\tau_y=\gamma2\delta1^\sim, \
\gamma\in\Cle{n-k-2},\ \delta\in\Cle{k-1}}
\wt(y)
+
\sum_{
y\in\mathcal Y(\nu), \
\tau_y=\alpha 12\xi1^\sim, \
\alpha\in\Cle{n-k-1},\ \xi\in\Cle{k-3}}
\wt(y).
\end{align*}

For a word~$y$ occurring in either sum, write~%
$\tau_y=t_1\dotsm t_r$.
Let~$t_1\dotsm t_q$ be the shortest prefix of~$\tau_y$
whose size is at least~$k-1$, namely, $u<k-1\le v$, where
$u=t_1+\dotsm+t_{q-1}$
and
$v=t_1+\dotsm+t_q$.
We claim that~$t_1\dotsm t_q\in\Cge{v}$. In fact, 
since~$I\in\mathcal C_{n-1}^{\le 2}$, 
every nonfinal part of~$\tau_y=I1^\sim$ is at least~$2$.
On the other hand, if~$t_q$ is the final part of~$\tau_y$, i.e., if~$q=r$, then 
\[
t_q=n-u\ge n-k+2\ge k+2\ge 5.
\]
This proves the claim.

Reading the column heights of~$I1$ from right to left gives~$\tau_y=(I1)^\sim$.
\Cref{fig:ribbon-crossing} illustrates the three cases (a), (b), and (c) below.
\begin{enumerate}[label=(\alph*)]
\item
In case~(a), writing~$\delta=2\delta'$ with~$\abs{\delta'}=k-3$ gives~%
$I1=\gamma22\delta'1$.
Then $v=1+\abs{\delta'}+1=k-1$,
and the next column has height~$t_{q+1}=2$,
coming from the two consecutive parts equal to~$2$ after~$\gamma$.
\item
In case~(b), writing~$\delta=1\delta'$ with~$\abs{\delta'}=k-2$ gives~%
$I1=\gamma21\delta'1$.
Then $v=2+\abs{\delta'}+1=k+1$,
and its last column has height~$t_q\ge 3$,
since it contains a box from each of the specified parts~$2$ and~$1$
and the first box of~$\delta'$.
\item
In case~(c), we have~$I1=\alpha12\xi1$ with~$\abs{\xi}=k-3$.
The indicated prefix has size $v=1+\abs{\xi}+1=k-1$,
and the next column has height~$t_{q+1}\ge 3$,
since it contains the last box of~$\alpha$
and a box from each of the specified parts~$1$ and~$2$.
\end{enumerate}
In each case, every row of~$I1$ has length at most~$2$,
and its last row has length~$1$.
Thus only the last run of~$y$ may have length~$1$, i.e.,~$y\in\mathcal Y$.
\begin{figure}[H]
\centering
\begin{tikzpicture}[x=.36cm,y=.36cm,font=\small,
hatched/.style={fill=black!12,
postaction={pattern=north east lines,pattern color=black}}]
\begin{scope}[xshift=0cm]
\node at (.3,7) {$\ddots$};
\node[anchor=east] at (-.1,6.6) {$\gamma$};
\draw (0,5) rectangle ++(1,1);
\draw[hatched] (1,5) rectangle ++(1,1);
\draw[hatched] (1,4) rectangle ++(1,1);
\draw[fill=black!12] (2,4) rectangle ++(1,1);
\node at (2.5,3) {$\ddots$};
\node at (3.5,1.8) {$\ddots$};
\draw[decorate,decoration={brace,amplitude=3pt}]
(5.3,4) -- (5.3,1)
node[midway,right=5pt] {$\delta'$};
\draw[fill=black!12] (4,0) rectangle ++(1,1);
\node[font=\scriptsize] at (4.5,.5) {$1$};
\draw[decorate,decoration={brace,mirror,amplitude=3pt}]
(2,-.6) -- (5,-.6);
\node at (2.5,-1.8) {(a)~$I1=\gamma22\delta'1$,~$\abs{\delta'}=k-3$};
\node at (2.5,-3.2) {$v=k-1,\quad t_{q+1}=2$};
\end{scope}
\begin{scope}[xshift=5.05cm]
\node at (.3,7) {$\ddots$};
\node[anchor=east] at (-.1,6.6) {$\gamma$};
\draw (0,5) rectangle ++(1,1);
\draw[fill=black!12] (1,5) rectangle ++(1,1);
\draw[fill=black!12] (1,4) rectangle ++(1,1);
\node[font=\scriptsize] at (1.5,3.5) {$\vdots$};
\draw[fill=black!12] (1,2) rectangle ++(1,1);
\node at (3,1.5) {$\ddots$};
\draw[decorate,decoration={brace,amplitude=3pt}]
(5.3,4) -- (5.3,1)
node[midway,right=5pt] {$\delta'$};
\draw[fill=black!12] (4,0) rectangle ++(1,1);
\node[font=\scriptsize] at (4.5,.5) {$1$};
\draw[decorate,decoration={brace,mirror,amplitude=3pt}]
(1,-.6) -- (5,-.6);
\node at (2.5,-1.8) {(b)~$I1=\gamma21\delta'1$,~$\abs{\delta'}=k-2$};
\node at (2.5,-3.2) {$v=k+1,\quad t_q\ge 3$};
\end{scope}
\begin{scope}[xshift=10.1cm]
\node at (.3,7.4) {$\ddots$};
\node[anchor=east] at (.4,6) {$\alpha$};
\draw[hatched] (1,6) rectangle ++(1,1);
\draw[hatched] (1,5) rectangle ++(1,1);
\draw[hatched] (1,4) rectangle ++(1,1);
\draw[fill=black!12] (2,4) rectangle ++(1,1);
\node at (2.5,3) {$\ddots$};
\node at (3.5,1.8) {$\ddots$};
\draw[decorate,decoration={brace,amplitude=3pt}]
(5.3,4) -- (5.3,1)
node[midway,right=5pt] {$\xi$};
\draw[fill=black!12] (4,0) rectangle ++(1,1);
\node[font=\scriptsize] at (4.5,.5) {$1$};
\draw[decorate,decoration={brace,mirror,amplitude=3pt}]
(2,-.6) -- (5,-.6);
\node at (2.5,-1.8) {(c)~$I1=\alpha12\xi1$,~$\abs{\xi}=k-3$};
\node at (2.5,-3.2) {$v=k-1,\quad t_{q+1}\ge 3$};
\end{scope}
\end{tikzpicture}
\caption{Ribbons appearing in~$\widetilde{U}_{n,k}$. 
Plain shaded boxes mark the prefix~$t_1\dotsm t_q$.
The hatched column has height~$t_{q+1}=2$ in (a) and~$t_{q+1}\ge 3$ in (c).
The side braces mark the ribbons~$\delta'$ in (a) and (b), and the ribbon~$\xi$ in (c).}
\label{fig:ribbon-crossing}
\end{figure}

Conversely, these three cases recover the corresponding compositions~$I$.
Thus
\begin{align*}
[s_\nu]\rho(\widetilde{U}_{n,k})
&=
\sum_{\theta\in\Cge{k-1}}
\norm{M_{\theta2}(\nu)}
+
\sum_{
\phi\in\Cge{k+1}, \
\last(\phi)\ge3}
\norm{M_\phi(\nu)}
+
\sum_{
\theta\in\Cge{k-1}, \
z\ge3}
\norm{M_{\theta z}(\nu)}
\\
&=
\sum_{\abs{\theta}=k-1}
\norm{M_\theta(\nu)}
+
\sum_{
\abs{\phi}=k+1, \
\last(\phi)\ge3}
\norm{M_\phi(\nu)}.
\end{align*}
In other words, this formula counts the total weights of words~$y\in\mathcal Y(\nu)$
such that~$y$ has a partial sum~$k-1$ or~$k+1$.

We now compute the contributions of the ribbon families indexed by~%
$\alpha22\beta1$,~$\alpha112\beta1$,~$\alpha21\beta1$, and~$\alpha111\beta1$.
Set~$c=2^{m_1(\alpha\beta)}$.
The four forms of~$I$ contain respectively~$0$,~$2$,~$1$, and~$3$
additional parts equal to~$1$ beyond those in~$\alpha\beta$, see
\cref{tab:ribbon-multipliers}.
Thus their coefficients~$2^{m_1(I)}$ in~$\widetilde{X}_{P_n}$ are~%
$c$,~$4c$,~$2c$, and~$8c$.
The two negative sums subtract~$2c$ from each of the first two forms
and~$3c$ from each of the last two forms.
Consequently,
\begin{align*}
\widetilde{E}_{n,k}-\widetilde{U}_{n,k}
&=
\smashoperator[r]{\sum_{
\alpha\in\Cle{n-k-1},\ \beta\in\Cle{k-4}}}
2^{m_1(\alpha\beta)}
\bigl(-R_{\alpha22\beta1}+2R_{\alpha112\beta1}\bigr)
+\smashoperator[r]{\sum_{
\alpha\in\Cle{n-k-1},\ \beta\in\Cle{k-3}}}
2^{m_1(\alpha\beta)}
\bigl(-R_{\alpha21\beta1}+5R_{\alpha111\beta1}\bigr).
\end{align*}
The final quotient divides the net coefficient by~$2^{m_1(I)}$,
because a coefficient~$2^{m_1(I)}$ of~$R_{I1}$ corresponds to
the weight~$\wt(y)$ of each corresponding Yamanouchi word.
The complete calculation is recorded in \cref{tab:ribbon-multipliers}.
\begin{center}
\begin{minipage}{\textwidth}
\centering
\small
\begin{tabular}{ccccc}
\toprule
Form of $I$
&
Positive coefficient
&
Coefficient subtracted
&
Net coefficient
&
Net multiplier
\\ \midrule
$\alpha22\beta$
&$c$&$2c$&$-c$&$(-c)/c=-1$
\\
$\alpha112\beta$
&$4c$&$2c$&$2c$&$(2c)/(4c)=1/2$
\\
$\alpha21\beta$
&$2c$&$3c$&$-c$&$(-c)/(2c)=-1/2$
\\
$\alpha111\beta$
&$8c$&$3c$&$5c$&$(5c)/(8c)=5/8$
\\ \bottomrule
\end{tabular}
\captionsetup{width=\textwidth}
\captionof{table}{Coefficients and net multipliers for the four ribbon families,
where~$c=2^{m_1(\alpha\beta)}$.}
\label{tab:ribbon-multipliers}
\end{minipage}
\end{center}

We next determine the restrictions on the run type
$\tau_y=(I1)^\sim=t_1\dotsm t_r$.
As before, set
$u=t_1+\dotsm+t_{q-1}$ and~$v=u+t_q$,
where~$t_q$ is the first part for which~$v\ge k-1$.
In each of the following figures, plain shading marks the columns~%
$t_1\dotsm t_q$, and the side brace marks the ribbon~$\beta$.
The four ribbon forms below are numbered (d)--(g).

\begin{enumerate}[label=(\alph*), start=4]
\item
\begin{minipage}[t]{.55\linewidth}
For~$I=\alpha22\beta$, the ribbon 
has~$\abs{\beta}=k-4$, see \cref{fig:ribbon-remainder:22}.
The columns strictly to the right of~$t_q$ contain all of~$\beta$,
the rightmost box of the second displayed part~$2$, and the final box.
Thus
\[
u=1+\abs{\beta}+1=k-2
\quad\text{and}\quad
v=k.
\]
Therefore the prefix~$t_1\dotsm t_q$ of~$\tau_y$ ends with a letter~$2$.
Since the net multiplier is~$-1$, this case contributes
\[
-\sum_{
\abs{\theta}=k, \
\last(\theta)=2}
\norm{M_\theta(\nu)}.
\]
\end{minipage}
\hfill
\begin{minipage}[t]{.43\linewidth}
\vspace{0pt}
\centering
\vspace{-\baselineskip}
\begin{tikzpicture}[x=.38cm,y=.38cm,font=\small]
\node at (.3,7) {$\ddots$};
\node[anchor=east] at (-.1,6.6) {$\alpha$};
\draw (0,5) rectangle ++(1,1);
\draw[fill=black!12] (1,5) rectangle ++(1,1);
\draw[fill=black!12] (1,4) rectangle ++(1,1);
\draw[fill=black!12] (2,4) rectangle ++(1,1);
\node at (3.2,2.4) {$\ddots$};
\draw[decorate,decoration={brace,amplitude=3pt}]
(6.5,3.8) -- (6.5,.9)
node[midway,right=5pt] {$\beta$};
\draw[fill=black!12] (5,0) rectangle ++(1,1);
\node[font=\scriptsize] at (5.5,.5) {$1$};
\end{tikzpicture}
\captionsetup{justification=centering}
\captionof{figure}{The ribbon~$I1=\alpha22\beta1$ in case~(d).}
\label{fig:ribbon-remainder:22}
\end{minipage}

\smallskip

The ribbons for cases~(e)--(g) are shown in
\cref{fig:ribbon-remainder:remaining}.

\item
For~$I=\alpha112\beta$, the ribbon again has~%
$\abs{\beta}=k-4$
and
$u=\abs{\beta}+2=k-2$.
The part~$t_q$ of~$\tau_y$ contains the two displayed parts~$1$ of~$I$,
the left box of the displayed part~$2$,
and at least the last box of~$\alpha$.
Hence~%
$t_q\ge4$ and~$v\ge k+2$.
Write~$z=t_q$.
With the net multiplier~$1/2$, this case contributes
\[
\frac{1}{2}
\sum_{
\abs{\theta}=k-2, \
z\ge4}
\norm{M_{\theta z}(\nu)}.
\]
\begin{figure}[H]
\vspace{-\baselineskip}
\centering
\begin{minipage}[t]{.31\linewidth}
\vspace{0pt}
\centering
\begin{tikzpicture}[x=.38cm,y=.38cm,font=\small]
\path[use as bounding box] (-.7,-.6) rectangle (7.5,8.2);
\node at (.4,7.8) {$\ddots$};
\node[anchor=east] at (.5,6.8) {$\alpha$};
\draw[fill=black!12] (1,6) rectangle ++(1,1);
\draw[fill=black!12] (1,5) rectangle ++(1,1);
\draw[fill=black!12] (1,4) rectangle ++(1,1);
\draw[fill=black!12] (1,3) rectangle ++(1,1);
\draw[fill=black!12] (2,3) rectangle ++(1,1);
\node at (3.3,2) {$\ddots$};
\draw[decorate,decoration={brace,amplitude=3pt}]
(6.5,3) -- (6.5,.9)
node[midway,right=5pt] {$\beta$};
\draw[fill=black!12] (5,0) rectangle ++(1,1);
\node[font=\scriptsize] at (5.5,.5) {$1$};
\end{tikzpicture}
\par\smallskip
\textup{(e)}~$I1=\alpha112\beta1$
\end{minipage}
\hfill
\begin{minipage}[t]{.31\linewidth}
\vspace{0pt}
\centering
\begin{tikzpicture}[x=.38cm,y=.38cm,font=\small]
\path[use as bounding box] (-.7,-.6) rectangle (7.5,8.2);
\node at (.3,7) {$\ddots$};
\node[anchor=east] at (-.1,6.6) {$\alpha$};
\draw (0,5) rectangle ++(1,1);
\draw[fill=black!12] (1,5) rectangle ++(1,1);
\draw[fill=black!12] (1,4) rectangle ++(1,1);
\fill (1.5,3.15) circle[radius=.4pt];
\fill (1.5,3.4) circle[radius=.4pt];
\fill (1.5,3.65) circle[radius=.4pt];
\draw[fill=black!12] (1,2) rectangle ++(1,1);
\draw[fill=black!12] (2,2) rectangle ++(1,1);
\node at (3.5,1.3) {$\ddots$};
\draw[decorate,decoration={brace,amplitude=3pt}]
(6.5,3.8) -- (6.5,.9)
node[midway,right=5pt] {$\beta$};
\draw[fill=black!12] (5,0) rectangle ++(1,1);
\node[font=\scriptsize] at (5.5,.5) {$1$};
\end{tikzpicture}
\par\smallskip
\textup{(f)}~$I1=\alpha21\beta1$
\end{minipage}
\hfill
\begin{minipage}[t]{.31\linewidth}
\vspace{0pt}
\centering
\begin{tikzpicture}[x=.38cm,y=.38cm,font=\small]
\path[use as bounding box] (-.7,-.6) rectangle (7.5,8.2);
\node at (.4,7.8) {$\ddots$};
\node[anchor=east] at (.5,6.8) {$\alpha$};
\draw[fill=black!12] (1,6) rectangle ++(1,1);
\draw[fill=black!12] (1,5) rectangle ++(1,1);
\draw[fill=black!12] (1,4) rectangle ++(1,1);
\draw[fill=black!12] (1,3) rectangle ++(1,1);
\fill (1.5,2.15) circle[radius=.4pt];
\fill (1.5,2.4) circle[radius=.4pt];
\fill (1.5,2.65) circle[radius=.4pt];
\draw[fill=black!12] (1,1) rectangle ++(1,1);
\draw[fill=black!12] (2,1) rectangle ++(1,1);
\node at (3.6,.8) {$\ddots$};
\draw[decorate,decoration={brace,amplitude=3pt}]
(6.5,3.8) -- (6.5,.9)
node[midway,right=5pt] {$\beta$};
\draw[fill=black!12] (5,-.5) rectangle ++(1,1);
\node[font=\scriptsize] at (5.5,0) {$1$};
\end{tikzpicture}
\par\smallskip
\textup{(g)}~$I1=\alpha111\beta1$
\end{minipage}
\caption{The ribbons for cases~(e)--(g).}
\label{fig:ribbon-remainder:remaining}
\end{figure}
\item
For~$I=\alpha21\beta$ with~$\abs{\beta}=k-3$, we have~%
$v=k$ and~$t_q\ge 3$; this case contributes
\[
-\frac{1}{2}
\sum_{
\abs{\theta}=k, \
\last(\theta)\ge3}
\norm{M_\theta(\nu)}.
\]
\item
For~$I=\alpha111\beta$ with~$\abs{\beta}=k-3$, 
we have~$u\le k-3$ and~$v\ge k+2$.
Set~$i=k-u$.
Then~$3\le i\le k$,
$t_1\dotsm t_{q-1}\in\Cge{k-i}$
and~%
$t_q\ge i+2$.
With net multiplier~$5/8$, this case contributes
\[
\frac{5}{8}
\sum_{i=3}^k
\sum_{
\abs{\theta}=k-i, \
z\ge i+2}
\norm{M_{\theta z}(\nu)}.
\]
\end{enumerate}
\Cref{tab:ribbon-run-cases} summarizes the seven classes and their multipliers.
\begin{table}[htbp]
\centering
\caption{Ribbon classes and their contributions per word of weight~$\wt(y)$.
Here~$u=t_1+\dots+t_{q-1}$ and~$v=u+t_q$, with~$q$ minimal such that~$v\ge k-1$.
All ribbon prefixes and suffixes have parts at most~$2$, and each ribbon has size~$n$.
Cases~(d) and~(e) are absent when~$k=3$.}
\label{tab:ribbon-run-cases}
\small
\begin{tabular}{cllc}
\toprule
Case & Ribbon~$I1$ and suffix size & Run condition & Multiplier \\
\midrule
(a) & $\gamma22\delta'1$,~$\abs{\delta'}=k-3$
& $v=k-1$,~$t_{q+1}=2$ & $1$ \\
(b) & $\gamma21\delta'1$,~$\abs{\delta'}=k-2$
& $v=k+1$,~$t_q\ge3$ & $1$ \\
(c) & $\alpha12\xi1$,~$\abs{\xi}=k-3$
& $v=k-1$,~$t_{q+1}\ge3$ & $1$ \\
(d) & $\alpha22\beta1$,~$\abs{\beta}=k-4$
& $u=k-2$,~$t_q=2$ & $-1$ \\
(e) & $\alpha112\beta1$,~$\abs{\beta}=k-4$
& $u=k-2$,~$t_q\ge4$ & $1/2$ \\
(f) & $\alpha21\beta1$,~$\abs{\beta}=k-3$
& $v=k$,~$t_q\ge3$ & $-1/2$ \\
(g) & $\alpha111\beta1$,~$\abs{\beta}=k-3$
& $u\le k-3$,~$v\ge k+2$ & $5/8$ \\
\bottomrule
\end{tabular}
\end{table}
Summing their contributions gives the desired formula.
\end{proof}

\section[Schur-positive spider families]{The spiders~\texorpdfstring{$S(a,2,1)$ and~$S(a,4,1)$}{S(a,2,1) and S(a,4,1)} are Schur positive}\label{sec:proof.s+}

In order to show that~$S(a, b, 1)$ is Schur positive, by \cref{X.Sab1:DE,cor:epos.A},
it suffices to prove the Schur positivity of the symmetric function~$E_{n,\,b+1}$.
We will show that~$E_{n, 3}$ and~$E_{n, 5}$ are Schur positive
in \cref{sec:Sa21,sec:Sa41} respectively.
For every positive integer~$n$, we write~$\underline{n}=12\dotsm n$
for the word obtained by concatenating~$1,2,\dots,n$ in that order.

\subsection[The family S(a,2,1)]{The spiders~\texorpdfstring{$S(a,2,1)$}{S(a,2,1)} are Schur positive}
\label{sec:Sa21}

By \cref{lem:s-coeff:E},
the Schur positivity of~$E_{n,3}$ follows from the~$j=3$ case
of the weighted run comparison in \cref{lem:multi-injection:j}.
For any word~$w$,
let~$\xi_k(w)$ denote the~$k$th run of~$w$.
We set~$\xi_k(w)=\emptyset$
if~$w$ has fewer than~$k$ runs.

\begin{lemma}\label{lem:multi-injection:j}
Let~$j\ge 3$ and~$n\ge j+1$.
Then for any partition~$\kappa\vdash n$,
\[
\norm{M_{j-1}(\kappa)}
\ge
\frac{1}{2}\norm{M_j(\kappa)}.
\]
\end{lemma}

\begin{proof}
We suppress~$\kappa$ from the notation in this proof.
Recall that~%
$\mathcal{Y}_\alpha$ is the set of Yamanouchi words~$y\in\mathcal{Y}$ 
whose run type~$\tau_y$ has prefix~$\alpha$.
Let~$y\in\mathcal{Y}_j$.
Then~$\xi_1(y)=\underline{j}$.
Indeed, the first run~$\xi_1(y)$ of~$y$ has length~$j$ by the definition of~$\mathcal{Y}_j$,
and the Yamanouchi inequalities force its strictly increasing letters to be~$1,2,\dots,j$.
We can therefore decompose
$\mathcal{Y}_j
=Y_1\sqcup Y_2\sqcup Y_3\sqcup Y_4$,
where
\begin{align*}
Y_1
&=\{y\in\mathcal{Y}_{j}\colon
j\not\in \xi_2(y)\},\\
Y_2
&=\{y\in\mathcal{Y}_{j}\colon
j\in \xi_2(y),\ 
\xi_3(y)\ne\emptyset,\ 
j\not\in \xi_3(y)\},\\
Y_3
&=\{y\in\mathcal{Y}_{j}\colon
j\in \xi_2(y),\ 
\xi_3(y)=\emptyset\},\quad\text{and}\\
Y_4
&=\{y\in\mathcal{Y}_{j}\colon
j\in \xi_2(y)\cap \xi_3(y)\}.
\end{align*}

We define a multi-map $f\colon M_j\to 2M_{j-1}$
as follows.
For~$y\in Y_1$,
let~$y_1$ be obtained by moving the terminal~$j$ of the first run
to the second run, and set~$f(y)=\{y_1\}$.
For~$y\in Y_2$,
let~$y_2$ be obtained by moving~$j$ from the first run
to the third run, and set~$f(y)=\{y_2\}$.
We use the symbol~$\mid$ to mark a boundary between consecutive runs. Then 
\[
Y_3=\begin{dcases*}
\{y_3\}, & if $n\ge 2j$ and $\kappa=(2^j,\,1^{n-2j})$, \\
\emptyset, & otherwise, 
\end{dcases*}
\]
where~$y_3=\underline{j}\mid \underline{n-j}$.
When~$Y_3$ is nonempty, define~$f(y_3)=\{y_{31}/2,\, y_{32}/2\}$,
where
\[
y_{31}=\underline{j-1}\mid1j\mid23\dotsm (n-j)
\quad\text{and}\quad
y_{32}=\underline{j-1}\mid\underline{j}\mid j(j+1)\dotsm (n-j).
\]
On the other hand, every word~$y_4\in Y_4$ has a unique factorization
\[
y_4=\underline{j}\mid1\alpha\mid1\beta\mid\zeta,
\]
where~$\alpha$ and~$\beta$ are increasing words beginning with~%
$23\dotsm j$,
and~$\zeta$ is the possibly empty suffix consisting of runs of~$y_4$.
Define $f(y_4)=\{y_{41}/2,\,y_{42}/2\}$,
where
\[
y_{41}=\underline{j-1}\mid1j\mid\alpha\mid1\beta\mid\zeta
\quad\text{and}\quad
y_{42}=\underline{j-1}\mid1j\mid1\alpha\mid\beta\mid\zeta.
\]
The definition of~$f$ is completed. 

First, we check that every target belongs to~$\mathcal Y_{j-1}$.
Each move preserves content.
For~$y_1$ and~$y_2$, only the count of~$j$ changes before its new position.
Delaying~$j$ preserves the inequality between the counts of~$j-1$ and~$j$.
For~$y_1$, increasing-order insertion places the moved~$j$ before any~$j+1$ in the second run;
no~$j+1$ occurs in the first run.
For~$y_2$, the second run begins with~$\underline j$;
after its~$j$ there is at most one~$j+1$ before the moved~$j$,
which is inserted before any~$j+1$ in the third run.
Thus the inequality between the counts of~$j$ and~$j+1$ also holds.

For~$y_{31}$ and~$y_{32}$, the displayed prefixes are Yamanouchi:
the letters~$1,\dots,j$ acquire their second occurrences in increasing order,
and any larger letters then occur for the first time in increasing order.
The prefix~$\underline{j-1}\mid1j\mid\alpha$ of~$y_{41}$
is Yamanouchi because~$\alpha$ begins with~$23\dotsm j$;
after this prefix, the counts agree with those in the first two runs of~$y_4$.
For~$y_{42}$, the extra~$1$ before~$\alpha$ only increases the count of~$1$;
after~$\underline{j-1}\mid1j\mid1\alpha$,
the counts agree with those in the first two runs of~$y_4$ followed by the initial~$1$ of its third run.
The remaining suffixes therefore satisfy the Yamanouchi inequalities as well.
The displayed boundaries remain boundaries between maximal increasing runs,
the first run has length~$j-1$, and every non-last run has length at least~$2$.
Hence all six targets lie in~$\mathcal Y_{j-1}$.

Next, the coefficients in every image sum to~$1$,
so~$f$ is mass-preserving.

Finally, we show that~$f$ is a multi-injection.
The four target forms $y_{31},y_{32},y_{41},y_{42}$ are pairwise distinct.
The words~$y_{31}$ and~$y_{32}$ have three runs,
whereas~$y_{41}$ and~$y_{42}$ have at least four.
The second runs of~$y_{31}$ and~$y_{32}$ are~$1j$ and~$\underline{j}$,
which are different since~$j\ge3$,
and the third runs of~$y_{41}$ and~$y_{42}$ are~$\alpha$ and~$1\alpha$,
which begin with~$2$ and~$1$, respectively.
Thus the run pattern determines which of the four forms a target has.

Moreover, none of $y_{31},y_{32},y_{41},y_{42}$ can arise as~$y_1$ or~$y_2$.
The only possible inverse for~$y_1$ moves~$j$ from the second run to the first,
whereas the only possible inverse for~$y_2$ moves~$j$ from the third run to the first.
Indeed,~$y_{31},y_{41},y_{42}$ all have second run~$1j$.
Applying the inverse operation for~$y_2$ to any of them
violates the Yamanouchi inequality between~$j-1$
and~$j$ at the end of the second run.
The inverse operation for~$y_1$ recovers a source in~$Y_3$ from~$y_{31}$,
recovers a source in~$Y_4$ from~$y_{41}$,
and produces a singleton non-last run from~$y_{42}$.
Thus none of these inverses lies in~$Y_1\sqcup Y_2$.
For~$y_{32}$,
either attempted inverse merges adjacent increasing factors and
recovers its source in~$Y_3$.

We can now count the preimages of each target.
The inverse operations above show that a target has at most one preimage
in each of~$Y_1$ and~$Y_2$.
For~$y_{31}$ and~$y_{32}$, the source is~$y_3$.
For~$y_{41}$ and~$y_{42}$, the run factorizations determine~%
$\alpha$,~$\beta$, and~$\zeta$, and hence~$y_4$, uniquely.
Together with the distinctness and separation just proved,
this shows that every target arising from~$Y_3\sqcup Y_4$
has exactly one preimage in~$\mathcal Y_j$ and occurs only once in its image.

These preimage counts give the required multiplicity bounds.
Let~$z\in\mathcal{Y}_{j-1}$.
If~$z$ has no preimage, then $m_{f(M_j)}(z)=0\le m_{2M_{j-1}}(z)$.
If~$z$ is the image of a source in~$Y_1$ or~$Y_2$,
then the defining move preserves the number of runs.
Hence every preimage has the same number of runs as~$z$ and the same weight.
Since~$z$ has at most one preimage in each of~$Y_1$ and~$Y_2$,
\[
m_{f(M_j)}(z)
\le
2\wt(z)
=m_{2M_{j-1}}(z).
\]
If~$z\in f(y)$ for some~$y\in Y_3\sqcup Y_4$,
then the displayed run factorizations show that~$z$ has one more run than~$y$,
so $\wt(y)=4\wt(z)$.
The uniqueness just proved and the coefficient~$1/2$ in~$f(y)$ give
\[
m_{f(M_j)}(z)
=\frac{1}{2}\wt(y)
=2\wt(z)
=m_{2M_{j-1}}(z).
\]
Therefore~$f$ is a mass-preserving multi-injection.

As a consequence, 
$2\norm{M_{j-1}}
=\norm{2M_{j-1}}
\ge
\norm{f(M_j)}
=\norm{M_j}$.
\end{proof}

Now let us show the first main result of this paper.

\begin{theorem}\label{thm:spos:Sa21}
Every spider of the form~$S(a,2,1)$ is Schur positive.
\end{theorem}

\begin{proof}
For~$a=1$, \cref{X.Sab1,thm:SW.path} give
\[
X_{S(1, 2, 1)}
=5e_5
+7e_{(4, 1)}
+e_{(3, 2)}
+2e_{(3, 1, 1)}
+e_{(2, 2, 1)},
\]
which is~$e$-positive and thus Schur positive.
Let~$a\ge 2$. Then the order~$n=a+4$ of~$S(a,2,1)$ is at least~$6$.
Let~$\nu\vdash n$,
and write~$M_\alpha=M_\alpha(\nu)$.
Set~$k=3$ in \cref{lem:s-coeff:E}.
Then~$n\ge 2k$,
\[
\Cge{0}=\{\emptyset\},
\quad
\Cge{1}=\emptyset,
\quad
\Cge{2}=\{2\},
\quad
\Cge{3}=\{3\},
\quad\text{and}\quad
\Cge{4}=\{4,\,22\}.
\]
By \cref{lem:multi-injection:j,lem:s-coeff:E},
\[
[s_\nu]E_{n,3}
\ge
\norm{M_{2}}
+\norm{M_{4}}
+\frac{5}{8}\sum_{z\ge 5}
\norm{M_{z}}
-\frac{1}{2}\norm{M_{3}}
\ge
\norm{M_{2}}
-\frac{1}{2}\norm{M_{3}}
\ge 0.
\]
Thus~$E_{n, 3}$ is Schur positive, and so is~$S(a, 2, 1)$.
\end{proof}

\subsection[The family S(a,4,1)]{The spiders~\texorpdfstring{$S(a,4,1)$}{S(a,4,1)} are Schur positive}\label{sec:Sa41}

As aforementioned, 
it suffices to prove that~$E_{n,5}$ is Schur positive for~$n\ge10$.
Our strategy is as follows.
We apply \cref{lem:s-coeff:E} with~$k=5$
and use content-preserving maps between Yamanouchi words
to reduce Schur positivity to inequalities between weighted multisets.
We verify these inequalities for words of length~$10$ by exhaustive computation,
keeping track of the final letter and run boundary.
We then extend them to longer words by appending common suffixes
and computing the resulting weight factors.

Let~$\mathcal{Y}(\kappa;\beta)$ denote the set of words in~$\mathcal{Y}(\kappa)$
having prefix~$\beta$, and set
\[
M(\kappa;\beta)=\brk1{\mathcal{Y}(\kappa;\beta),\,\wt}.
\]
Let~$M_\alpha(\kappa;\beta)$ be the induced sub-multiset of~$M_\alpha(\kappa)$
on words having prefix~$\beta$.
For compactness,
set 
\[
\omega=121233
\quad\text{and}\quad
\omega_z=\underline{3}\mid \underline{z}.
\]

\begin{lemma}\label{lem:k=5}
Let~$n\ge 6$ and~$\kappa\vdash n$.
The following statements hold.
\begin{enumerate}[label=(\arabic*)]
\item
There is a mass-preserving multi-injection
\[
\iota\colon M_{32}(\kappa)\to M_{23}(\kappa)
\]
such that
$\wt(\iota(y))=\wt(y)$
and
$M_{23}(\kappa)\backslash
\iota(M_{32}(\kappa))
=M(\kappa;\omega)$.
\item
For any~$z\ge 4$, if~$n\ge z+2$, then
there is a mass-preserving multi-injection
\[
\iota\colon
M_{2z}(\kappa)
\to M_{3(z-1)}(\kappa)
\]
such that
$\wt(\iota(y))=\wt(y)$
and
$
M_{3(z-1)}(\kappa)\backslash\iota(M_{2z}(\kappa))
=M_{3(z-1)}(\kappa;\omega_{z-1})$.
\item
If~$n\ge 7$, then
there is an injection of underlying sets
\[
\varphi\colon
\underline{M(\kappa;\omega)}\to
\bigsqcup_{z\ge 4}\underline{M_{3z}(\kappa;\omega_z)}
\]
such that
$\wt(\varphi(y))=4\wt(y)$
for all~$y\in M(\kappa;\omega)$.
As a consequence, 
\[
\sum_{z\ge4}\norm{M_{3z}(\kappa;\omega_z)}
\ge 4\norm{M(\kappa;\omega)}.
\]
\end{enumerate}
\end{lemma}

\begin{proof}
Fix~$\kappa$ and suppress it from the notation.
In the formulas below,~$\zeta$ is a possibly empty suffix.

\begin{enumerate}[label=(\arabic*)]
\item
For~$y\in\mathcal Y_{32}$, the first two runs satisfy~%
$\xi_1=123$ and~$\xi_2\in\{12,14\}$.
Define
\[
\iota(123\mid12\mid\zeta)=12\mid123\mid\zeta
\quad\text{and}\quad
\iota(123\mid14\mid\zeta)=12\mid134\mid\zeta.
\]
Then~$\iota$ preserves the weight by \cref{wt.y}.
Conversely, every word in~$\mathcal Y_{23}$ begins with~%
$12\mid123$ or~$12\mid134$.
Thus the inverse of~$\iota$ is well defined 
except when a word starts with~$\omega=121233$.
Thus~$\iota$ is a weight-preserving bijection onto~%
$\mathcal Y_{23}\backslash\mathcal Y(\kappa;\omega)$,
which gives the required mass-preserving multi-injection and complement identity.

\item
For~$y\in\mathcal Y_{2z}$, the first two runs satisfy~%
$\xi_1=12$ and $\xi_2\in\{\underline z,\,134\dotsm(z+1)\}$.
Define
\[
\iota(12\mid\underline z\mid\zeta)
=123\mid124\dotsm z\mid\zeta
\quad\text{and}\quad
\iota(12\mid134\dotsm(z+1)\mid\zeta)
=123\mid14\dotsm(z+1)\mid\zeta.
\]
Both replacements preserve the final letter of the second run,
all subsequent run boundaries, and the weight.
The Yamanouchi inequalities force the second run of a word in~%
$\mathcal Y_{3(z-1)}$ to be one of~$\underline{z-1}$,~$124\dotsm z$,
and~$14\dotsm(z+1)$.
The last two forms have unique preimages under the displayed replacements.
The first gives exactly the words with prefix
$\omega_{z-1}=123\mid\underline{z-1}$.
This proves the mass-preserving multi-injection and complement identity.

\item
Since~$n\ge7$, every word with prefix~$\omega$ has a unique run factorization
\[
y=12\mid123\mid3\dotsm z\mid\zeta
\quad\text{with }z\ge4.
\]
Define $\varphi(y)=123\mid\underline z\mid\zeta$.
Then $\varphi(y)\in M_{3z}(\,;\omega_z)$,
since the unchanged final letter~$z$ preserves the boundary with~$\zeta$.
The second run of~$\varphi(y)$ determines~$z$,
and the displayed replacement uniquely recovers~$y$.
Thus~$\varphi$ is injective.
It preserves length and decreases the number of runs by one,
so $\wt(\varphi(y))=4\wt(y)$.
Summing over the image gives the stated norm inequality.
\qedhere
\end{enumerate}
\end{proof}

Let~$n\ge 10$ and~$\kappa\vdash n$.
Define two multisets
\[
X(\kappa)
=\brk1{
\mathcal Y_{23}(\kappa)
\backslash\mathcal Y(\kappa;\omega),\,
3\wt/2}
\quad\text{and}\quad
Y(\kappa)
=\brk1{2M_{24}(\kappa)}
\sqcup
M(\kappa;1231231)
\sqcup
M_{22}(\kappa).
\]
For~$b\in\{<,\,\ge\}$ and~$t\in\mathbb Z_{\ge 1}$,
let~$X_b(\kappa,t)$ be the induced sub-multiset of~$X(\kappa)$
on words satisfying~$y_9\mathrel b y_{10}$ and~$y_{10}=t$, 
where~$y_9$ and~$y_{10}$ mean the ninth and tenth letter of~$y$.
Define~$Y_b(\kappa,t)$ analogously.
These restrictions retain the multiplicities in~$X(\kappa)$ and~$Y(\kappa)$.

\begin{proposition}\label{prop:10}
For any content~$\mu\vdash 10$ and any letter~$t\in \mathbb Z_{\ge 1}$,
\begin{align*}
\norm{X_<(\mu,t)}\le \norm{Y_<(\mu,t)}
\quad\text{and}\quad
\norm{X_\ge(\mu,t)}\le \norm{Y_\ge(\mu,t)}.
\end{align*}
\end{proposition}

\begin{proof}
The exhaustive verification using exact integer arithmetic appears in \cref{sec:finite-verification}.
\end{proof}

We remark that the length\nobreakdash-$10$ threshold is best possible for these inequalities.
For this length\nobreakdash-$9$ example only,
define~$X(\kappa)$ and~$Y(\kappa)$ by the same formulas,
and let~$X_b(\kappa,t)$ and~$Y_b(\kappa,t)$ denote their induced sub-multisets
on words satisfying~$y_8\mathrel b y_9$ and~$y_9=t$.
Then, for~$\kappa=432$ and~$t=2$,
\[
X_<(\kappa,t)
=\{121231312\}
\quad\text{and}\quad
Y_<(\kappa,t)
=\{121312312\}.
\]
Hence
$\norm{X_<(\kappa,t)}=6>4=\norm{Y_<(\kappa,t)}$.

\begin{lemma}[Finite-prefix extension]\label{lem:finite-prefix}
For every~$n\ge10$ and every partition~$\nu\vdash n$,
\[
\norm{X(\nu)}\le\norm{Y(\nu)}.
\]
\end{lemma}

\begin{proof}
For~$n=10$, this follows by summing \cref{prop:10} over~$t$ and~$b$.
Suppose that~$n\ge11$.
Fix~$\mu\vdash10$,~$1\le t\le10$,~$b\in\{<,\,\ge\}$,
and a word~$\zeta$ of length~$n-10$ and componentwise content~$\nu-\mu$.
We compare the words~$w\zeta$ obtained by taking~$w$ from~%
$X_b(\mu,t)$ or~$Y_b(\mu,t)$.

The ballot inequalities for~$w\zeta$ beyond~$w$ are
\[
\mu_i+m_i(\zeta')\ge\mu_{i+1}+m_{i+1}(\zeta')
\qquad(i\ge1)
\]
for each prefix~$\zeta'$ of~$\zeta$, with missing parts of~$\mu$ taken as zero.
They depend only on~$\mu$ and~$\zeta$.
For the run-length condition, compare~$t$ with the first letter~$\zeta_1$ of~$\zeta$.
If~$t<\zeta_1$, the last run of~$w$ and the first run of~$\zeta$
merge into one run of length at least~$2$.
If~$t\ge\zeta_1$, the last run of~$w$ becomes a non-last run of~$w\zeta$,
so it must have length at least~$2$; this holds precisely when~$b$ is~$<$.
The other non-last runs in~$w$ are already admissible,
and those in~$\zeta$ are fixed.
Thus either all these concatenations lie in~$\mathcal Y(\nu)$ or none do.

In the admissible case, appending~$\zeta$ preserves the first-seven-letter conditions defining~$X$ and~$Y$
and their weight multipliers~$3/2$,~$2$, and~$1$.
By \cref{wt.y} and the two run cases above,
\[
\frac{\wt(w\zeta)}{\wt(w)}
=
\begin{cases}
2^{n-8-2\ell(\tau_\zeta)},&t<\zeta_1,\\
2^{n-10-2\ell(\tau_\zeta)},&t\ge\zeta_1.
\end{cases}
\]
For fixed~$\mu$,~$t$,~$b$, and~$\zeta$, this positive factor is independent of~$w$.
Multiplying the corresponding inequality in \cref{prop:10} by this factor bounds the total multiplicity of the extensions from~$X_b(\mu,t)$
by that of those from~$Y_b(\mu,t)$.
Every word in~$X(\nu)$ or~$Y(\nu)$ has a unique such decomposition:
its length\nobreakdash-$10$ prefix satisfies the ballot and run conditions
and belongs to the corresponding~$X_b(\mu,t)$ or~$Y_b(\mu,t)$.
Summing over~$\mu$,~$t$,~$b$, and~$\zeta$ proves the claim.
\end{proof}

\begin{theorem}\label{thm:spos:Sa41}
For every integer~$a\ge1$,
the spider~$S(a,4,1)$ is Schur positive.
\end{theorem}

\begin{proof}
For~$a=1$, the graph~$S(1,4,1)\cong S(4,1,1)$ is Schur positive by \citet{WW23-DAM}.
For~$a=2$, the graph~$S(2,4,1)\cong S(4,2,1)$ is covered by \cref{thm:spos:Sa21}.
For~$a=3$, the graph~$S(3,4,1)\cong S(4,3,1)$ is~$e$-positive by \citet{WW23-DAM},
and therefore Schur positive.
Now let~$a\ge4$, and let~$n=a+6\ge10$ be the number of vertices of~$S(a,4,1)$.
Let~$\nu\vdash n$.
By \cref{X.Sab1:DE},
it suffices to prove that~$[s_\nu]E_{n,5}\ge 0$.
For the moment,
we suppress the parameter~$\nu$.
By \cref{lem:s-coeff:E},
\begin{align*}
[s_\nu]E_{n,5}
&\ge
\brk2{\norm{M_{4}}
+\norm{M_{22}}}
+\brk2{\norm{M_{6}}
+\norm{M_{33}}
+\norm{M_{24}}}
+\frac{1}{2}\sum_{z\ge 4}
\norm{M_{3z}}\\
&\qquad
+\frac{5}{8}
\brk3{
\sum_{z\ge 5}
\norm{M_{2z}}
+\sum_{z\ge 7}
\norm{M_{z}}}
-\norm{M_{32}}
-\frac{1}{2}
\brk2{\norm{M_{5}}
+\norm{M_{23}}}.
\end{align*}
Since $\norm{M_4}\ge \norm{M_5}/2$ by \cref{lem:multi-injection:j},
and since~$\norm{M_z}\ge 0$,
we obtain
\begin{equation}\creflabel[ineq]{ineq:Sa41-reduction}
[s_\nu]E_{n,5}
\ge
\norm{M_{22}}
+\brk2{\norm{M_{33}}
+\norm{M_{24}}}
+\frac{1}{2}\sum_{z\ge 4}
\norm{M_{3z}}
-\brk2{\norm{M_{32}}
+\frac{1}{2}\norm{M_{23}}}.
\end{equation}

By \cref{lem:k=5}(2) with~$z=4$,
\begin{equation}\label{eq:Sa41-M33-M24}
\norm{M_{33}}
=\norm{M_{24}}+\norm{M_{33}(\nu;\omega_3)}
=\norm{M_{24}}+\norm{M(\nu;1231231)}.
\end{equation}
Applying \cref{lem:k=5}(2) with~$z+1$ in place of~$z$ gives
\[
\norm{M_{3z}}
\ge \norm{M_{3z}(\nu;\omega_z)}.
\]
Also, \cref{lem:k=5}(1) gives
$\norm{M_{32}}=\norm{M_{23}\backslash M(\nu;\omega)}$, and
\begin{equation}\label{eq:Sa41-M32-M23}
\norm{M_{32}}+\frac12\norm{M_{23}}
=\frac32\norm{M_{23}\backslash M(\nu;\omega)}
+\frac12\norm{M(\nu;\omega)}.
\end{equation}
Finally, \cref{lem:k=5}(3) gives
\begin{equation}\creflabel[ineq]{ineq:Sa41-prefix-bound}
\frac12\sum_{z\ge4}\norm{M_{3z}(\nu;\omega_z)}
-\frac12\norm{M(\nu;\omega)}
\ge\frac32\norm{M(\nu;\omega)}
\ge 0.
\end{equation}
Substituting \cref{eq:Sa41-M33-M24,eq:Sa41-M32-M23}
into \cref{ineq:Sa41-reduction}
and applying \cref{ineq:Sa41-prefix-bound} yields
\begin{align*}
[s_\nu]E_{n,5}
&\ge
\norm{M_{22}}
+2\norm{M_{24}}
+\norm{M(\nu;1231231)}
-\frac{3}{2}
\norm{M_{23}\backslash M(\nu;\omega)}\\
&=
\norm{Y(\nu)}-\norm{X(\nu)}.
\end{align*}
By \cref{lem:finite-prefix}, the last difference is nonnegative, which proves the theorem.
\end{proof}

\appendix
\crefalias{section}{appendix}
\section{A sign-reversing involution proof of the path expansion}\label{sec:path-involution}

We give another proof of the coefficient formula in \cref{thm:path.ribbon},
using hook decompositions.
Let $I=i_1\dotsm i_s$ and~$J$ be compositions of~$n$.
The \emph{ribbon decomposition of~$J$ relative to~$I$} is a decomposition
\[
J=J_1\bullet\dotsm\bullet J_s
\]
in which each~$J_r$ is a ribbon,
the sequence of their sizes is the composition~$I$,
and each binary operation~$\bullet$ is either the concatenation or the near concatenation.
A \emph{hook} is a ribbon of the form~$1^s t$ for some~$s\ge 0$ and~$t\ge 1$.
A \emph{hook decomposition} is a ribbon decomposition into hooks.

\citet[Proposition~4.23]{GKLLRT95} gives the following transition formula.

\begin{proposition}\label{prop:hook-transition}
For any composition~$I=i_1\dotsm i_s$, we have
\[
\varepsilon^I\Psi^I
=\sum_{J=J_1\bullet\dotsm\bullet J_s}
\varepsilon^{J_1\dotsm J_s}R_J,
\]
where the sum runs over all hook decompositions of~$J$ relative to~$I$.
\end{proposition}

\begin{proof}[Alternative proof of the coefficient formula in \cref{thm:path.ribbon}]
Using the power-sum expansion established in the proof of \cref{thm:SW.path}
and the hook transition formula in \cref{prop:hook-transition},
\[
\widetilde{X}_{P_n}
=\sum_{I\vDash n}\varepsilon^I \Psi^I
=\sum_{I\vDash n}
\sum_{J=J_1\bullet J_2\bullet\dotsm}
\varepsilon^{J_1J_2\dotsm}
R_J
=\sum_{J\vDash n}
\sum_{J_1\bullet J_2\bullet\dotsm\in\mathcal{H}(J)}
\varepsilon^{J_1 J_2\dotsm}
R_J, 
\]
where the inner sum in the middle expression runs over all hook decompositions of~$J$ relative to~$I$,
and~$\mathcal{H}(J)$ is the set of all hook decompositions of~$J$.
We will simplify this expression by means of a sign-reversing involution.

Let~$J\vDash n$.
A \emph{row tail} of~$J$ is the first box met in a row of length at least~$2$ when the ribbon is traversed from its northeastern end toward its southwestern end.
A \emph{leader} of a hook decomposition is a box that begins a hook of length at least~$2$.
Let~$\mathcal{H}'(J)$ be the set of hook decompositions of~$J$ in which every row tail is a leader.
We construct an involution~$\varphi$ on~$\mathcal{H}(J)\backslash \mathcal{H}'(J)$.

Let $d\in\mathcal{H}(J)\backslash \mathcal{H}'(J)$.
Let~$R$ be the first row encountered in this northeast-to-southwest traversal whose tail is not a leader in~$d$.
Let~$\Box$ and~$\Box'$ be the first two boxes of~$R$ in that order.
Define~$\varphi(d)$ to be the hook decomposition obtained from~$d$ by the following local replacement.
If they lie in one hook,
that hook is necessarily of the form~$1^{s-1}(t+1)$ for some~%
$s,t\ge1$;
replace it locally by the two hooks~$1^s\triangleright t$.
If they lie in distinct hooks,
those hooks are necessarily~$1^s$ and~$t$, joined by~$\triangleright$;
replace~$1^s\triangleright t$ by the single hook~%
$1^{s-1}(t+1)$.
These two local configurations are inverse to one another.
They leave every earlier row unchanged and leave~$R$ defective:
if~$t=1$, its tail begins only a singleton hook,
whereas if~$t>1$, its tail does not begin a hook.
Thus the same row is selected after either operation and~%
$\varphi^2(d)=d$.
Moreover,
the involution~$\varphi$ is sign reversing since
\[
\varepsilon^{1^{s-1}(t+1)}
=(-1)^t
=-\varepsilon^{1^s}\varepsilon^t.
\]

If~$d\in\mathcal H'(J)$,
no hook can have a horizontal arm:
the terminal box of such an arm would be a row tail that is not a leader.
Consequently all hooks are vertical,
which forces every row of~$J$ to have length~$1$ or~$2$ and its last row
to have length~$1$.
These are exactly the defining conditions for~$J\in\Cle{n-1}1$.
Conversely,
if~$J\in\Cle{n-1}1$,
the members of~$\mathcal H'(J)$ are precisely the decompositions into
vertical hooks in which the first two boxes of every non-leftmost column
remain in one hook.
It follows that the coefficient~$[R_J]\widetilde{X}_{P_n}$ vanishes unless~$J\in\Cle{n-1}1$.

For~$J\in\Cle{n-1}1$,
every hook~$J_i$ in every decomposition in~$\mathcal{H}'(J)$ is vertical,
and therefore has sign~$\varepsilon^{J_i}=1$.
Hence
\begin{equation}\label{pf:XPn:H'}
\widetilde{X}_{P_n}
=\sum_{J\in\Cle{n-1}1}
\sum_{J_1\bullet J_2\bullet\dotsm\in\mathcal{H}'(J)}
\varepsilon^{J_1}\varepsilon^{J_2}\dotsm
R_J
=\sum_{J\in\Cle{n-1}1}
\abs{\mathcal{H}'(J)}
R_J.
\end{equation}
Let~$l=\ell(J^\sim)$ be the number of columns.
A vertical-hook decomposition chooses cuts among the~$n-l$ vertical adjacencies.
Each of the~$l-1$ non-leftmost columns has height at least~$2$
and forbids a cut between its first two boxes.
The remaining~$n-2l+1$ cuts are independent, so \cref{lem:m1:ell.sim} gives
\[
\abs{\mathcal H'(J)}
=2^{n-2l+1}
=2^{m_1(J)-1}.
\]
Substituting this into \cref{pf:XPn:H'},
we obtain
\[
\widetilde{X}_{P_n}
=\frac{1}{2}\sum_{J\in\Cle{n-1}1}
2^{m_1(J)}R_J
=\sum_{I\in\Cle{n-1}}
2^{m_1(I)}R_{I1}.\qedhere
\]
\end{proof}

\section{Exact finite verification}\label{sec:finite-verification}

We verify \cref{prop:10} by generating Yamanouchi words recursively from the empty word.
For a prefix with counts~$c_1\ge\dots\ge c_h>0$, set~$c_{h+1}=0$
and append each letter~$i\in\{1,\dots,h+1\}$ such that~$i=1$ or~$c_{i-1}>c_i$.
This is exactly the ballot condition for the extended prefix,
so the recursion generates every Yamanouchi word once.
At length~$10$, it gives~$9{,}496$ words, of which~$1{,}366$ have no non-last run of length~$1$.
For each retained word~$y$ with~$r$ runs, \cref{wt.y} gives~$\wt(y)=2^{11-2r}$.
Since~$10\ge2(r-1)+1$, we have~$r\le5$;
the weight is even, so all multiplicities in~$X$ and~$Y$ are exact integers.

Summing these multiplicities verifies all~$840$ inequalities
($42$ contents,~$10$ final letters, and two relations):
there are~$86$ strict inequalities and~$754$ equalities with both sides zero.
For~$t>10$, both sides are zero.
In \cref{tab:finite-comparisons}, write~$X_b=\norm{X_b(\mu,t)}$ and~$Y_b=\norm{Y_b(\mu,t)}$.
The table gives all~$54$ pairs with~$X_<+X_\ge>0$;
both left sides vanish for every omitted pair.
Partition exponents indicate multiplicities.

To reproduce the verification with Python~$3.9$ or later and its standard library, run
\begin{center}
\texttt{python3 -B }\path{scripts/verify_finite_comparison.py}
\end{center}
The accompanying~\path{scripts/README.md} describes the complete CSV certificate,
recorded output, checksums, and table regeneration.

\begingroup
\floatsetup[table]{capposition=top}
\begin{table}[H]
\centering
\small
\setlength{\tabcolsep}{4pt}
\caption{All pairs with a nonzero left side in the finite comparison.
For every listed pair, $X_<\le Y_<$ and $X_\ge\le Y_\ge$.}\label{tab:finite-comparisons}
\begin{tabular}{@{}crrrrr@{\hspace{4em}}crrrrr@{\hspace{4em}}crrrrr@{}}
\toprule
$\mu$ & $t$ & $X_<$ & $Y_<$ & $X_\ge$ & $Y_\ge$ &
$\mu$ & $t$ & $X_<$ & $Y_<$ & $X_\ge$ & $Y_\ge$ &
$\mu$ & $t$ & $X_<$ & $Y_<$ & $X_\ge$ & $Y_\ge$ \\
\midrule
$541$ & 1 & 0 & 0 & 3 & 4 & $431^{3}$ & 5 & 48 & 84 & 0 & 0 & $3^{2}1^{4}$ & 2 & 12 & 24 & 30 & 56 \\
$532$ & 1 & 0 & 0 & 6 & 8 & $42^{3}$ & 1 & 0 & 0 & 18 & 26 & $3^{2}1^{4}$ & 6 & 108 & 196 & 0 & 0 \\
$531^{2}$ & 1 & 0 & 0 & 9 & 12 & $42^{3}$ & 4 & 48 & 68 & 3 & 6 & $32^{3}1$ & 1 & 0 & 0 & 15 & 38 \\
$52^{2}1$ & 1 & 0 & 0 & 9 & 12 & $42^{2}1^{2}$ & 1 & 0 & 0 & 36 & 66 & $32^{3}1$ & 4 & 60 & 158 & 30 & 46 \\
$521^{3}$ & 1 & 0 & 0 & 9 & 12 & $42^{2}1^{2}$ & 3 & 36 & 62 & 6 & 10 & $32^{3}1$ & 5 & 96 & 142 & 0 & 0 \\
$51^{5}$ & 1 & 0 & 0 & 3 & 4 & $42^{2}1^{2}$ & 5 & 60 & 92 & 0 & 0 & $32^{2}1^{3}$ & 1 & 0 & 0 & 15 & 44 \\
$4^{2}2$ & 2 & 12 & 14 & 6 & 8 & $421^{4}$ & 1 & 0 & 0 & 30 & 56 & $32^{2}1^{3}$ & 3 & 24 & 48 & 30 & 56 \\
$4^{2}2$ & 3 & 12 & 14 & 0 & 2 & $421^{4}$ & 2 & 12 & 24 & 3 & 4 & $32^{2}1^{3}$ & 6 & 180 & 342 & 0 & 0 \\
$4^{2}1^{2}$ & 2 & 12 & 30 & 9 & 12 & $421^{4}$ & 6 & 60 & 102 & 0 & 0 & $321^{5}$ & 1 & 0 & 0 & 12 & 24 \\
$4^{2}1^{2}$ & 4 & 12 & 22 & 0 & 0 & $41^{6}$ & 1 & 0 & 0 & 12 & 24 & $321^{5}$ & 2 & 0 & 0 & 12 & 24 \\
$43^{2}$ & 1 & 0 & 0 & 3 & 12 & $41^{6}$ & 7 & 24 & 40 & 0 & 0 & $321^{5}$ & 7 & 144 & 272 & 0 & 0 \\
$43^{2}$ & 3 & 24 & 26 & 3 & 6 & $3^{3}1$ & 3 & 36 & 48 & 21 & 40 & $31^{7}$ & 8 & 48 & 96 & 0 & 0 \\
$4321$ & 1 & 0 & 0 & 24 & 42 & $3^{3}1$ & 4 & 24 & 52 & 0 & 0 & $2^{5}$ & 5 & 24 & 98 & 12 & 20 \\
$4321$ & 2 & 24 & 40 & 9 & 12 & $3^{2}2^{2}$ & 2 & 12 & 18 & 18 & 26 & $2^{4}1^{2}$ & 4 & 12 & 24 & 15 & 44 \\
$4321$ & 3 & 36 & 52 & 3 & 8 & $3^{2}2^{2}$ & 4 & 72 & 148 & 15 & 18 & $2^{4}1^{2}$ & 6 & 72 & 178 & 0 & 0 \\
$4321$ & 4 & 36 & 52 & 0 & 0 & $3^{2}21^{2}$ & 2 & 12 & 42 & 36 & 66 & $2^{3}1^{4}$ & 3 & 0 & 0 & 12 & 24 \\
$431^{3}$ & 1 & 0 & 0 & 21 & 36 & $3^{2}21^{2}$ & 3 & 36 & 60 & 18 & 34 & $2^{3}1^{4}$ & 7 & 72 & 200 & 0 & 0 \\
$431^{3}$ & 2 & 24 & 52 & 9 & 12 & $3^{2}21^{2}$ & 5 & 132 & 220 & 0 & 0 & $2^{2}1^{6}$ & 8 & 48 & 96 & 0 & 0 \\
\bottomrule
\end{tabular}
\end{table}
\endgroup


\bibliographystyle{spos-Thibon-bibstyle}
\IfFileExists{csf.bib}
{\bibliography{csf}}
{\IfFileExists{../csf.bib}
{\bibliography{../csf}}
{\bibliography{csf}}}

\end{document}